\documentclass[12pt, a4paper]{article}
\usepackage{amsfonts}
\usepackage{mathrsfs}
\usepackage{latexsym}
\usepackage{xy}
\usepackage{amsfonts,amsmath,amssymb,amsthm}
\usepackage{color}

\xyoption{all}
\newcommand{\mput}{\multiput}
\newcommand{\bcen}{\begin{center}}     \newcommand{\ecen}{\end{center}}
\newcommand{\bay}{\begin{array}}      \newcommand{\eay}{\end{array}}
\newcommand{\beq}{\begin{eqnarray*}}      \newcommand{\eeq}{\end{eqnarray*}}
\def\az{\alpha}
\def\bz{\beta}

\def\i{\mathscr{I}}
\def\s{\mathscr{S}}

\def\ch{\mathrm{char}}

\def\gr{\mathrm{gr}}
\def\Gr{\mathrm{Gr}}
\def\Hom{\mathrm{Hom}}

\def\op{\mathrm{op}}
\def\ot{\otimes}
\def\Ext{\mathrm{Ext}}

\def\End{\mathrm{End}}

\def\dim{\mathrm{dim}}
\def\mod{\mathrm{mod}}
\def\Mod{\mathrm{Mod}}

\def\Im{\mathrm{Im}}
\def\Ker{\mathrm{Ker}}

\def\Rep{\mathrm{Rep}}
\def\rep{\mathrm{rep}}

\def\top{\mathrm{top}}

\def\z{\mathbb{Z}}
\begin{document}

\newtheorem{theorem}{Theorem}[section]
\newtheorem{proposition}{Proposition}[section]
\newtheorem{lemma}{Lemma}[section]
\newtheorem{corollary}{Corollary}[section]
\newtheorem{remark}{Remark}[section]
\newtheorem{example}{Example}[section]
\newtheorem{definition}{Definition}[section]

\title{Superspecies and their representations
\footnote{Project 10731070 supported by NSFC.}}

\author{Yang Han $^1$ and Deke Zhao $^2$}

\date{\footnotesize 1. Academy of Mathematics and Systems Science,
Chinese Academy of Sciences,\\ Beijing 100080, P.R.China. E-mail:
hany@iss.ac.cn\\ 2. School of Mathematical Science, Beijing Normal
University, Beijing 100875,\\ P.R. China. E-mail:
dkzhao@mail.bnu.edu.cn}

\maketitle

\begin{abstract} Superspecies are introduced to provide
the nice constructions of all finite-dimensional superalgebras. All
acyclic superspecies, or equivalently all finite-dimensional
(gr-basic) gr-hereditary superalgebras, are classified according to
their graded representation types. To this end, graded equivalence,
graded representation type and graded species are introduced for
finite group graded algebras.
\end{abstract}

{\footnotesize {\bf Mathematics Subject Classification 2000}: 16W50,
16G20, 16W55, 16G60}

\section{Introduction}

\indent\indent Clifford algebras, more special Grassmann algebras,
and their representations play a quite important role in both
mathematics and physics (ref. \cite{V} and some references therein).
Here, we shall study more general finite-dimensional superalgebras
and their representation theory.

It is well-known that quiver and species play a crucial role in the
construction and representation theory of finite-dimensional
algebras (ref. \cite{ARS,Ga2,DR2}). Is there an analog of quiver and
species for finite-dimensional superalgebras? In order to answer
this question, we shall introduce graded equivalence theory at first
(Section 3). Our aim is to reduce the constructions and the
representation theory of all finite-dimensional finite group graded
algebras to those of all gr-basic ones. We shall show that every
finite-dimensional finite group graded algebra is graded equivalent
to a gr-basic one, which is unique up to congruence
(Theorem~\ref{treducedtobasic}). Furthermore, we shall introduce the
concept of graded species (Section 4.2), in particular superspecies,
which plays a crucial role for finite-dimensional superalgebras as
quiver and species for finite-dimensional algebras. In case the
underlying field is algebraically closed and of characteristic not
equal to 2, we shall obtain the super version of Wedderburn's
principal theorem (Theorem~\ref{tWedderburn}), and show that each
finite-dimensional gr-basic superalgebra is isomorphic to the factor
of the graded tensor algebra of a superspecies modulo an admissible
graded ideal (Theorem~\ref{tpresentation}). In particular, each
finite-dimensional gr-basic gr-hereditary superalgebra is graded
isomorphic to the graded tensor algebra of an acyclic superspecies.

Next we consider the question:``How to classify acyclic
superspecies, or equivalently finite-dimensional (gr-basic)
gr-hereditary superalgebras, according to their graded
representation types?'' In order to define the graded representation
types of acyclic superspecies, we shall introduce the graded
representation types of finite-dimensional finite group graded
algebras at first (Section 5). Then we shall provide the graded
version of Drozd's theorem (Theorem~\ref{tDrozd}). After classifying
all finite-dimensional gr-division superalgebras over an
algebraically closed field, we shall introduce the quiver of a
superspecies (Section 6.1), and prove that the category of
finite-dimensional representations of a superspecies is equivalent
to the category of finite-dimensional representations of its quiver
(Theorem~\ref{tspeciesquiver}). Furthermore, we shall obtain the
super version of Kac's theorem (Theorem~\ref{tKactheorem}) and
classify all acyclic superspecies according to their graded
representation types in terms of their quivers
(Theorem~\ref{tclassification1}). For a superspecies, we shall also
introduce its superquiver (Section 6.3). The construction of the
superquiver of a superspecies is easier than that of its quiver.
Therefore, we shall classify all acyclic superspecies again
according to their graded representation types in terms of their
superquivers (Theorem~\ref{tclassification2}).

In fact, our original motivation is to realize some Lie
superalgebras and their quantized enveloping superalgebras by
finite-dimensional gr-hereditary superalgebras, just as the
realization of Kac-Moody algebras and their quantized enveloping
algebras by finite-dimensional hereditary algebras via Hall algebra
approach (ref. \cite{R,Gr,PX}). As far as we know, many experts
thought or are thinking about this problem. But we think, before
this, we should learn much more about finite-dimensional
superalgebras and their representation theory.

\section{Graded algebras and graded modules}

In this section, on one hand, we shall fix some notations and
terminologies on graded algebras and graded modules, on the other
hand, we shall get some new results on these aspects as well. For
the knowledge of graded ring theory, we refer to \cite{NV}.

\subsection{Graded algebras}

\indent\indent Throughout we assume that $K$ is a fixed field, $G$
is a finite multiplicative group with identity $e$, and the
composition of maps is written from left to right.

A {\it $G$-graded $K$-vector space} or {\it graded vector space} is
a $K$-vector space $V = \oplus _{g \in G} V_g$ where $V_g$'s are
$K$-subspaces of $V$. By definition, each $v \in V$ can be uniquely
written as $v=\sum_{g \in G}v_g$ where $v_g \in V_g$ for all $g \in
G$. A nonzero element $v \in V_g$ is said to be {\it homogeneous of
degree $g$}, and we write $\deg v =g$. A graded vector space $V$ is
said to be {\it trivially graded} if $V = V_e$. A {\it graded
subspace} of a graded vector space $V$ is a $K$-subspace $U$ of $V$
such that $U = \oplus _{g \in G} (U \cap V_g)$. A {\it map} of
graded vector space is a $K$-linear map $\phi : V \rightarrow W$
such that $\phi(V_g) \subseteq W_g$ for all $g \in G$. If $G$ is
abelian, it is usually considered as an additive group. In
particular, in the case of $G = \mathbb{Z}_2 = \mathbb{Z}/2
\mathbb{Z}$, a $G$-graded $K$-vector space is called a {\it super
vector space}.

A {\it $G$-graded $K$-algebra} or {\it graded algebra} $A$ is both
an associative $K$-algebra with identity and a $G$-graded $K$-vector
space $A = \oplus_{g \in G} A_g$ such that $A_gA_h \subseteq A_{gh}$
for all $g, h \in G$. A {\it graded subalgebra} (resp. {\it graded
ideal}) of a graded algebra is both a subalgebra (resp. ideal) and a
graded subspace. A {\it graded homomorphism} of graded algebras is
both an algebra homomorphism and a graded vector space map. In the
case of $G = \mathbb{Z}_2$, a $G$-graded $K$-algebra is called a
{\it superalgebra}.

\subsection{Graded modules}

\indent\indent A (left) {\it $G$-graded module} or {\it graded
module} $M$ over a graded algebra $A = \oplus _{g \in G} A_g$ is
both a (left) $A$-module and a graded vector space $M = \oplus _{g
\in G} M_g$ such that $A_gM_h \subseteq M_{gh}$ for all $g, h \in
G$. A {\it graded morphism} of graded modules is both an $A$-module
morphism and a graded vector space map.

All graded $A$-modules and all graded morphisms between them form a
Grothendieck category, denoted by $\Gr A$. Denote by $\gr A$ the
full subcategory of $\Gr A$ consisting of all finitely generated
graded $A$-modules. Moreover, we denote by $\Mod A$ the category of
left (ungraded) $A$-modules, and by $\mod A$ the full subcategory of
$\Mod A$ consisting of finitely generated left $A$-modules. Denote
by ${\cal F}$ the {\it forgetful functor} from $\gr A$ (resp. $\Gr
A$) to $\mod A$ (resp. $\Mod A$).

The {\it $g$-th shift functor} or {\it $g$-suspension functor} $S_g$
from $\Gr A$ or $\gr A$ to itself, is a functor defined on objects
by $S_g(X)=Y$ with $Y_h=X_{hg}$, and in the obvious way on
morphisms. Clearly, $S_g$ is an automorphism functor and
$S_gS_h=S_{gh}$ for all $g,h \in G$.

Let $M , N \in \Gr A$. Define $\Hom_A(M,N)_g := \{\phi \in
\Hom_{\Mod A}({\cal F}(M),{\cal F}(N)) \mid \phi(M_h) \subseteq
N_{hg}, \forall h \in G \}$ and $\Hom_A(M,N):=\Hom_{\Mod A}({\cal
F}(M),{\cal F}(N))$. Clearly $\Hom_A(M,N)_g = {\cal F}\Hom_{\Gr
A}(M,S_g(N)) = {\cal F}\Hom_{\Gr A}(S_{g^{-1}}(M),N)$. Non-zero
elements of $\Hom_A(M,N)_g$ are called {\it graded $A$-morphisms of
degree $g$}. Since $G$ is a finite group, we have
$\Hom_A(M,N)=\oplus_{g \in G}\Hom_A(M,N)_g$ (ref. \cite[Corollary
3.10]{Da}). Furthermore, $\End_A(M) := \Hom_A(M,M)$ is a graded
algebra. Denote by $\Ext_{\Gr A}^n(-, N)$ (resp. $\Ext_A^n(-, N)$)
the $n$-th right derived functor of the functor $\Hom_{\Gr A}(-, N)$
(resp. $\Hom_A(-, N)$). Since $G$ is a finite group, we have
$\Ext^n_A(M,N)=\oplus_{g \in G}\Ext^n_A(M,N)_g$ where
$\Ext^n_A(M,N)_g \cong \Ext^n_{\Gr A}(M,S_g(N))$ as $K$-vector
spaces for all $g \in G$.

\subsection{Gr-indecomposable modules}

\indent\indent Let $A$ be a graded algebra. A nonzero graded
$A$-module $M$ is said to be {\it gr-indecomposable} if it is not
the direct sum of two nonzero graded modules.

Let $KG^*$ be the dual vector space of $KG$ and $\{p_g|g \in G\}$
its dual basis such that $p_g(x)=a_g \in K$ for all $g \in G$ and
$x=\sum_{h \in G}a_hh \in KG$. The {\it smash product} $A \# KG^*$
is the $K$-algebra whose underlying $K$-vector space is $A \otimes
_K KG^*$ and whose multiplication is given by $(a \# p_g)(b \# p_h)
:= ab_{gh^{-1}} \# p_h$ where $a \# p_g$ denotes $a \otimes p_g$
(ref. \cite[Section 1]{CM}). It is well-known that $\Gr A \cong \Mod
A \# KG^*$ (ref. \cite[Theorem 2.2]{CM}). By this isomorphism and
the Krull-Schmidt theorem for finite-dimensional algebras, or using
the same strategy as the proof of \cite[Corollary 5.3]{P}, we obtain
the following result:

\begin{theorem} \label{tkrullschmidt} {\bf (Graded version of Krull-Schmidt theorem)}
Let $A$ be a finite-dimensional graded algebra and $M \in \gr A$.
Then $M$ is gr-indecomposable if and only if $\End_{\gr A}(M)$ is
local. Moreover, each nonzero graded module in $\gr A$ has a unique
decomposition of gr-indecomposables up to permutation and graded
isomorphism. \end{theorem}

A graded module $M$ is said to be {\it gr-simple} if $0$ and $M$ are
the only graded submodules of $M$. A graded module $M$ is said to be
{\it gr-semisimple} if $M$ is a direct sum of gr-simple modules. A
graded submodule $N$ of a graded module $M$ is called a {\it
gr-maximal submodule} if $M/N$ is gr-simple. For a graded module
$M$, we denote by $J_G(M)$ its {\it graded Jacobson radical}, i.e.,
the intersection of all gr-maximal submodules of $M$. We call
$\top_G(M):=M/J_G(M)$ the {\it gr-top} of $M$. We say that a graded
submodule $N$ of a graded module $M$ is {\it gr-small} in $M$ if
$N+X=M$ for a graded submodule $X$ of $M$ implies $X=M$. Clearly,
$N$ is gr-small in $M$ if and only if $N \subseteq J_G(M)$. Note
that $J_G({_A}A)$ is a graded ideal of $A$, called the {\it graded
Jacobson radical} of $A$, and denoted by $J_G(A)$. Moreover, $J_G(A)
\subseteq J(A)$, the Jacobson radical of $A$ (ref. \cite[Theorem
4.4(1)]{CM}).

The following result is a graded analog of \cite[Proposition
I.3.5]{ARS}, which can be proved similarly:

\begin{proposition}\label{pradical} Let $A$ be a finite-dimensional
graded algebra and $M \in \gr A$. Then $J_G(M)=J_G(A)M$.
\end{proposition}

A graded algebra is called a {\it gr-division algebra} if every
nonzero homogeneous element is invertible. A graded algebra is said
to be {\it gr-local} if every homogeneous element is either
invertible or nilpotent.

The following result provides a characterization of gr-local
algebras:

\begin{proposition} \label{pgradedlocal} Let $A$ be a finite-dimensional
graded algebra. Then the following statements are equivalent:

(1) $A$ is gr-local.

(2) $A/J_G(A)$ is a gr-division algebra.

(3) The initial subalgebra $A_e$ of $A$ is local. \end{proposition}

\begin{proof} $(1) \Rightarrow (2)$: Suppose that $A$ is gr-local.
Let $J$ be the graded ideal of $A$ generated by all
homogeneous nilpotent elements in $A$.

{\bf Claim 1.} If $x \in A_g$ is nilpotent then $axb$ is nilpotent
for all homogeneous elements $a, b$ in $A$.

{\it Proof of Claim 1.} First of all, we show that $ax$ is
nilpotent. Assume on the contrary that $ax$ is invertible. Since $A$
is graded local, $a$ is either invertible or nilpotent. If $a$ is
invertible then so is $x=a^{-1}(ax)$. It is a contradiction. If $a$
is nilpotent then there is $n \geq 2$ such that $a^n=0$ and $a^{n-1}
\neq 0$. On one hand, we have $a^nx=0$. On the other hand, since
$ax$ is invertible, we have $a^nx = a^{n-1}(ax) \neq 0$. It is also
a contradiction. Similarly, one can show that $xb$ is nilpotent.
Thus $axb$ is nilpotent as well.

{\bf Claim 2.} If $x, y \in A_g$ are nilpotent then $x+y$ is
nilpotent.

{\it Proof of Claim 2.} If $x=0$ then we need do nothing. If $x \neq
0$ then assume on the contrary that $x+y$ is invertible, i.e., there
is $z \in A_{g^{-1}}$ such that $(x+y)z=1$. Since $x$ is nilpotent,
there is $n \geq 2$ such that $x^n=0$ and $x^{n-1} \neq 0$. On one
hand, we have $x^{n-1}(1-yz)=0$. On the other hand, by Claim 1, $yz$
is nilpotent, thus $1-yz$ is invertible and $x^{n-1}(1-yz) \neq 0$.
It is a contradiction.

It follows from Claim 1 and Claim 2 that all homogeneous elements in
$J$ are nilpotent. Thus $A/J$ is a gr-division algebra and $J$ is
the unique gr-maximal ideal of $A$. Hence $J = J_G(A)$ and
$A/J_G(A)$ is a gr-division algebra.

$(2) \Rightarrow (1)$: Let $x$ be a homogeneous element of $A$. We
show that $x$ is either nilpotent or invertible. Since $J_G(A)
\subseteq J(A)$, $J_G(A)$ is nilpotent. If $x \in J_G(A)$ then $x$
is nilpotent. If $x \notin J_G(A)$, since $A/J_G(A)$ is a
gr-division algebra, $\bar{x}$ is invertible in $A/J_G(A)$. So is
$x$ by \cite[Proposition 2.9.1 vi)]{NV}.

$(1) \Rightarrow (3)$: It is trivial.

$(3)\Rightarrow (1)$: Suppose that $A_e$ is local. For any
homogeneous element $x$ of $A$, $x^{|G|} \in A_e$ is either
invertible or nilpotent. So is $x$.
\end{proof}

\begin{remark} The conclusion $(1) \Leftrightarrow (3)$ is similar
to \cite[Theorem 3.1]{GG}, which says that a $\mathbb{Z}$-graded
Artin algebra with local initial subring is local. However, the
latter does not hold for finite group graded algebras: We consider
the $\mathbb{Z}_2$-graded algebra $A = A_{\bar{0}} \oplus
A_{\bar{1}}$ where $A_{\bar{0}} := K, A_{\bar{1}} := K \varepsilon,$
and $\varepsilon^2 := 1$. Clearly, the initial subalgebra of $A$ is
local, but $A$ is not local in the case of $\ch K \neq 2$, since
$\frac{1}{2}(1+\varepsilon)$ is an idempotent which is neither
invertible nor nilpotent.

This example also implies that \cite[Theorem 3.2]{GG}, which says
that for a $\mathbb{Z}$-graded Artin algebra $A$, $M \in \gr A$ is
gr-indecomposable if and only if ${\cal F}(M)$ is an indecomposable
$A$-module, does not hold for finite group graded algebras: Indeed,
in the case of $\ch K \neq 2$, $A$ is itself a gr-simple module but
${\cal F}(A)$ is decomposable, since forgotten grading $A \cong K
\oplus K$ as algebras. \end{remark}

Now we give another characterization of gr-indecomposable modules.

\begin{proposition} \label{pindecomposable} Let $A$ be a finite-dimensional graded
algebra and $M \in \gr A$. Then $M$ is gr-indecomposable if and only
if $\End_A(M)$ is gr-local. Moreover, if $M \in \gr A$ is gr-simple
then $\End_A(M)$ is a gr-division algebra.
\end{proposition}

\begin{proof} By Theorem~\ref{tkrullschmidt} and
Proposition~\ref{pgradedlocal}, $M$ is gr-indecomposable if and only
if $\End_A(M)$ is gr-local.

Let $M$ be gr-simple. For any $g \in G$ and $\phi = {\cal F}(\psi)
\in \End_A(M)_g$ with $\psi \in \Hom_{\Gr A}(M,S_g(M))$, since $M$
and $S_g(M)$ are gr-simple, we have $\Im \psi = 0$ or $S_g(M)$, and
$\Ker \psi = 0$ or $M$. If $\Im \psi = 0$ or $\Ker \psi = M$ then
$\psi = 0$. Otherwise, $\Im \psi = S_g(M)$ and $\Ker \psi = 0$,
i.e., $\psi$ and thus $\phi$ is invertible.
\end{proof}

\subsection{Gr-projective modules}

\indent\indent We say $F \in \Gr A$ is {\it gr-free} if it has an
$A$-basis consisting of homogeneous elements, equivalently $F \cong
\oplus_{i \in I}S_{g_i}(A)$ with $g_i \in G$ for all $i \in I$. A
graded module $P$ is said to be {\it gr-projective} if for any
graded epimorphism $\phi: M \rightarrow N$ and any graded morphism
$\psi: P \rightarrow N$ there is a graded morphism $\varphi: P
\rightarrow M$ such that $\psi=\varphi\phi$.

\begin{proposition} \label{pprojective} Let $A$ be a graded algebra and $P \in \Gr A$. Then the
following assertions are equivalent:

(1) $P$ is gr-projective;

(2) $\Hom_{\Gr A}(P, -): \Gr A \rightarrow \Mod K$ is exact;

(3) $\Hom_A(P, -): \Gr A \rightarrow \Mod K$ is exact;

(4) $P$ is a direct summand of a gr-free module;

(5) $\Ext^1_{\Gr A}(P, N)=0$  for all $N \in \Gr A$;

(6) $\Ext^1_A(P, N)=0$ for all  $N \in \Gr A$.
\end{proposition}

\begin{proof} $(1) \Leftrightarrow (2)$: By the definition of gr-projectiveness.

$(2) \Rightarrow (4)$: Applying $\Hom_{\Gr A}(P, -)$ to a graded
epimorphism $\phi : F \rightarrow P$, we obtain $\phi$ splits.

$(4) \Rightarrow (3)$: Since ${\cal F}(P)$ is projective, the
functor $\Hom_{\Mod A}({\cal F}(P), -): \Mod A \rightarrow \Mod K$
is exact. Thus the functor $\Hom_A(P, -): \Gr A \rightarrow \Mod K$
is exact.

$(3) \Rightarrow (6)$: By the definition of extension group, we are
done.

$(6) \Rightarrow (5)$: It follows from $\Ext^n_A(M, N) = \oplus_{g
\in G}\Ext^n_A(M, N)_g, \forall M, N \in \Gr A$.

$(5) \Leftrightarrow (2)$: It is clear by $\Gr A \cong \Mod A \#
KG^*$. \end{proof}

Now we classify all finite-dimensional gr-indecomposable
gr-projective modules over a finite-dimensional graded algebra.

\begin{proposition} \label{pgradedprojective} Let $A$ be a finite-dimensional
graded algebra. Then, up to graded isomorphism, all
finite-dimensional gr-indecomposable gr-projective $A$-modules are
of form $S_g(P)$ where $g \in G$ and $P$ is a gr-indecomposable
direct summand of $A$. \end{proposition}

\begin{proof} Let $M$ be a finite-dimensional gr-indecomposable
gr-projective $A$-modules. Since $M$ is finitely generated, $M =
\sum_{i=1}^n A m_i$ for some homogeneous elements $m_i$ of degree
$g_i$. Thus $\phi : \oplus_{i=1}^n S_{g_i^{-1}}(A) \rightarrow M,
(a_1,...,a_n) \mapsto \sum_{i=1}^n a_im_i$ is a graded $A$-module
epimorphism. Since $M$ is gr-projective, $M$ is graded isomorphic to
a direct summand of $\oplus_{i=1}^n S_{g_i^{-1}}(A)$. By
Theorem~\ref{tkrullschmidt}, we are done. \end{proof}

Next we classify all finite-dimensional gr-simple modules over a
finite-dimensional graded algebra.

\begin{proposition} \label{pgradedsimple} Let $A$ be a finite-dimensional
graded algebra. Then, up to graded isomorphism, all
finite-dimensional gr-simple $A$-modules are of form
$S_g(P)/J_G(S_g(P)) \cong S_g(P/J_G(P))$ where $g \in G$ and $P$ is
a \linebreak gr-indecomposable direct summand of $A$.
\end{proposition}

\begin{proof} First of all, we show that for any finite-dimensional
gr-indecomposable gr-projective $A$-module $P$, $P/J_G(P)$ is
gr-simple. Indeed, there is a natural algebra epimorphism $\End_{\Gr
A}(P) \rightarrow \End_{\Gr A}(P/J_G(P))$. Since $P$ is
gr-indecomposable, $\End_{\Gr A}(P)$ is local. So is $\End_{\Gr
A}(P/J_G(P))$. By Theorem~\ref{tkrullschmidt}, $P/J_G(P)$ is
gr-indecomposable. Since $P/J_G(P)$ is gr-semisimple, it must be
gr-simple.

Conversely, for any finite-dimensional gr-simple $A$-module $M$,
there is a finite-dimensional gr-projective $A$-module $P$ and a
graded $A$-module epimorphism $\phi : P \rightarrow M$. Since $M$ is
gr-simple, $J_G(P) \subseteq \Ker \phi$. Thus there is a graded
$A$-module epimorphism $\bar{\phi} : P/J_G(P) \rightarrow M$.
Decompose $P$ into the direct sum of gr-indecomposable gr-projective
modules, say $P= \oplus_{i=1}^nP_i$. Then $P/J_G(P) \cong
\oplus_{i=1}^nP_i/J_G(P_i)$ and there is a graded $A$-module
epimorphism $\psi : \oplus_{i=1}^nP_i/J_G(P_i) \rightarrow M$.
Assume that the restriction of $\psi$ on $P_i/J_G(P_i)$ is nonzero.
Since both $P_i/J_G(P_i)$ and $M$ are gr-simple, we have
$P_i/J_G(P_i) \cong M$.
\end{proof}

\subsection{Gr-hereditary algebras}

\indent\indent Let $A$ be a finite-dimensional graded algebra. The
\textit{gr-projective dimension} of $M \in \Gr A$, denoted by
$\mbox{gr-pd}_A M$, is the least integer $n$ for which there exists
an exact sequence $P_n \rightarrowtail \cdots \rightarrow P_0
\twoheadrightarrow M$ in $\Gr A$ where $P_0, \cdots, P_n$ are
gr-projective. If no such sequence exists for any $n$ then
$\mbox{gr-pd}_A M := \infty$ (cf. \cite[\S2.3]{NV}). The
\textit{gr-global dimension} of $A$ is $\mbox{gr-gl.dim} A := \sup
\{\mbox{gr-pd}_A M | M \in \Gr A\}$ (cf. \cite[\S6.3]{NV}).

\begin{proposition} \label{pprojdim} Let $A$ be a finite-dimensional graded
algebra and $M \in \Gr A$. Then the following assertions are
equivalent:

(1) $\mbox{gr-pd}_AM \leq n$;

(2) $\Ext_{\Gr A}^{n+1}(M, N)=0$ for all $ N \in \Gr A$;

(3) $\Ext_A^{n+1}(M, N)=0$ for all $N \in \Gr A$.
\end{proposition}

\begin{proof}
$(1) \Rightarrow (3):$ It follows from $\mbox{gr-pd}_AM \le n$ that
$\mbox{pd}_AM \le n$. Thus $\Ext_{\Mod A}^{n+1}({\cal F}(M), {\cal
F}(N))=0$ and further $\Ext_A^{n+1}(M, N)=0$ for all $N \in \Gr A$.

$(3) \Rightarrow (2):$ It follows from $\Ext_A^{n+1}(M, N) =
\oplus_{g \in G}\Ext_{\Gr A}^{n+1}(M, N)_g$.

$(1) \Leftrightarrow (2):$ It is clear by $\Gr A \cong \Mod A \#
KG^*$. \end{proof}

\begin{corollary} \label{cgldim} Let $A$ be a finite-dimensional graded
algebra. Then the following assertions are equivalent:

(1) $\mbox{gr-gl.dim}A \le n$;

(2) $\Ext_{\Gr A}^{n+1}(M,N)=0$ for all $M, N \in \Gr A$;

(3) $\Ext_A^{n+1}(M,N)=0$ for all $M, N \in \Gr A$.
\end{corollary}

The same strategy as the proof of \cite[Theorem 1]{Au} can be used
to characterize gr-global dimension as follows:

\begin{proposition}\label{pgldim} Let $A$ be a finite-dimensional graded
algebra. Then

$\begin{array}{rcl} \mbox{gr-gl.dim}A
&=&\sup \{\mbox{gr-pd}_A A/I| I \text{ is a graded left ideal of } A \}\\
&=&\mbox{gr-pd}_A (A/J_G(A))\\ & = & \sup \{\mbox{gr-pd}_AS| S
\text{ is gr-simple}\}.\end{array}$
\end{proposition}

\begin{proposition} \label{phereditary}
Let $A$ be a finite-dimensional graded algebra. Then the following
conditions are equivalent:

(1) All graded left ideals of $A$ are gr-projective;

(2) $J_G(A)$ is gr-projective;

(3) $\mbox{gr-pd} _A(A/J_G(A))\le 1$;

(4) $\mbox{gr-gl.dim} A \le 1$.
\end{proposition}

\begin{proof} $(1) \Rightarrow (2) \Rightarrow(3)$: Trivial.

$(3)\Rightarrow(4):$ It follows from Proposition~\ref{pgldim}.

$(4)\Rightarrow(1):$ For any graded ideal $I$ of $A$, it follows
from Corollary~\ref{cgldim} that $\Ext_{\Gr A}^2(A/I,N)=0, \forall N
\in \Gr A$. Thus $\Ext_{\Gr A}^1(I,N)=0, \forall N \in \Gr A$. By
Proposition~\ref{pprojective}, we are done.
\end{proof}

\begin{definition} We say a finite-dimensional graded algebra is
\mbox{\bf gr-hereditary} if it satisfies the equivalent conditions
of Proposition~\ref{phereditary}. \end{definition}

A graded epimorphism $\varphi: M \rightarrow N$ is said to be {\it
gr-essential} if a graded morphism $\psi: X \rightarrow M$ is a
graded epimorphism whenever $\varphi\psi: X \rightarrow N$ is a
graded epimorphism.

\begin{proposition} \label{pessential} Let $A$ be a finite-dimensional graded algebra and $\varphi:
M \rightarrow N$ a graded epimorphism in $\gr A$. Then the following
assertions are equivalent:

(1) $\varphi$ is gr-essential;

(2) $\Ker \varphi \subseteq J_G(M)$;

(3) The induced graded epimorphism $\bar{\varphi}:
M/J_G(M)\rightarrow N/J_G(N)$ is a graded
isomorphism.\end{proposition}

\begin{proof} (1) $\Rightarrow$ (2): Suppose that $\varphi$ is
gr-essential and $X$ is a graded submodule of $M$ such that $X+ \Ker
\varphi=M$. Denote by $\lambda: X\hookrightarrow M$ the natural
graded monomorphism. Then  $\varphi\lambda: X\rightarrow N$ is a
graded epimorphism. Thus $\lambda$ is a graded epimorphism. Hence
$X=M$ and $\Ker \varphi$ is gr-small. So $\Ker \varphi \subseteq
J_G(M)$.

(2) $\Rightarrow$ (3): Since $\varphi$ is a graded epimorphism, one
has a natural graded isomorphism $\tilde{\varphi}: M/\Ker \varphi
\rightarrow N$. It follows from $\Ker \varphi \subseteq J_G(M)$ that
$J_G(M/\Ker \varphi)= J_G(M)/\Ker \varphi$. Thus $\tilde{\varphi}$
induces a graded isomorphism $\bar{\varphi}: M/J_G(M)=(M/\Ker
\varphi)/(J_G(M)/\Ker \varphi) \rightarrow N/J_G(N)$.

(3) $\Rightarrow$ (1): Suppose $\psi: X\rightarrow M$ is a graded
morphism such that $\varphi\psi: X\rightarrow N$ is a graded
epimorphism. Then $\bar{\varphi}\bar{\psi}: X/J_G(X)\rightarrow
N/J_G(N)$ is a graded epimorphism. Since $\bar{\varphi}$ is a graded
isomorphism, $\bar{\psi}$ is a graded epimorphism. Thus $\Im \psi
+J_G(M)=M$. Hence $\Im \psi=M$, i.e., $\psi$ is a graded
epimorphism.
\end{proof}

A {\it gr-projective cover} of a graded module $M$ is a gr-essential
graded epimorphism $\varphi: P \rightarrow M$ with $P$
gr-projective. The following result can be obtained from $\gr A
\cong \mod A \# KG^*$ and \cite[Theorem I.4.2]{ARS}:

\begin{proposition} \label{pprojcover} Let $A$ be a finite-dimensional
graded algebra. Then each $M \in \gr A$ has a gr-projective cover
which is unique up to isomorphism, namely, if $\varphi_i : P_i
\rightarrow M, i=1,2,$ are two gr-projective covers of $M$ then
there is a graded isomorphism $\phi : P_1 \rightarrow P_2$ such that
$\varphi_2\phi=\varphi_1$.
\end{proposition}

The following result is an immediate consequence of
Proposition~\ref{pessential}:

\begin{corollary}\label{cprojcovertop} A graded epimorphism $\varphi: P \rightarrow M$ with $P$
gr-projective is a gr-projective cover if and only if the induced
graded epimorphism $\top_G(P) \rightarrow \top_G(M)$ is a graded
isomorphism.
\end{corollary}

\section{Graded equivalence}

In this section, we shall introduce graded equivalence theory so
that we can reduce the constructions and the representation theory
of finite-dimensional finite group graded algebras to those of
gr-basic ones. Note that Gordon and Green had introduced the graded
equivalence theory for $\mathbb{Z}$-graded Artin algebras (ref.
\cite[Section 5]{GG}). For finite-dimensional finite group graded
algebras, it is very similar. However, the introduction of graded
natural transformation is new and necessary, which seems to be
neglected before.

\subsection{Graded equivalences}

\indent\indent Let $A$ and $B$ be graded algebras. A functor ${\cal
U}: \Gr A \rightarrow \Gr B$ is called a {\it graded functor} if it
commutes with the shift functors $S_g$ for all $g \in G$.

Let $M , N \in \Gr A$. Then $\End_A(M)$ is a $G$-graded algebra and
$\Hom_A(M,N)$ is a graded $\End_A(M)$-modules. For any $\az \in
\Hom_A(N,L)_e$, $\Hom_A(M,\az): \Hom_A(M,N) \rightarrow
\Hom_A(M,L)$, $(\bz_g)_g \mapsto (\bz_g \az)_g$, is a graded
$\End_A(M)$-morphism of degree $e$. Thus we obtain a functor
$\Hom_A(M,-): \Gr A \rightarrow \Gr(\End_A(M))$.

\begin{lemma} \label{lgradehom} Let $A$ be a $G$-graded algebra
and $M \in \Gr A$. Then $\Hom_A(M,-): \Gr A \rightarrow \Gr
(\End_A(M))$ is a graded functor such that ${\cal F}\Hom_A(M,-) =$
\linebreak $\Hom_{\Mod A}({\cal F}(M),-){\cal F}$.
\end{lemma}

Let $A$ and $B$ be graded algebras and ${\cal U}, {\cal V} : \Gr A
\rightarrow \Gr B$ two graded functors. A natural transformation
$\eta : {\cal U} \rightarrow {\cal V}$ is said to be {\it graded} if
$\eta_{S_g(M)} = S_g(\eta_M)$ for all $g \in G$ and $M \in \Gr A$.
We say that a graded functor ${\cal U}: \Gr A \rightarrow \Gr B$ is
a {\it graded equivalence}, if there is a graded functor ${\cal V}:
\Gr B \rightarrow \Gr A$ such that there are graded natural
isomorphisms ${\cal V}{\cal U} \cong 1_{\Gr A}$ and ${\cal U}{\cal
V} \cong 1_{\Gr B}$. We say that $A$ and $B$ are {\it graded
equivalent} if there is a graded equivalence $\Gr A \rightarrow \Gr
B$.

Let ${\cal L}: \Mod A \rightarrow \Mod B$ be an equivalence. We call
${\cal L}$ a {\it graded equivalence} if there is a graded functor
${\cal U}: \Gr A \rightarrow \Gr B$ such that ${\cal F}{\cal U} =
{\cal L}{\cal F}$. In this case ${\cal U}$ is called an {\it
associated graded functor} of ${\cal L}$.

\begin{proposition} \label{pmoritaequ} Let ${\cal U}$ be an associated graded functor of
a graded equivalence ${\cal L}: \Mod A \rightarrow \Mod B$, $Q={\cal
U}(A)$, and $P=\Hom_B(Q, B)$. Then, as graded algebras, $A \cong
\End_B(Q)$ and $B \cong \End_A(P)$. Moreover, the functor ${\cal
M}:=\Hom_{\Mod B}({\cal F}(Q), -)$ is inverse to ${\cal L}$, ${\cal
L} \cong {\cal L}' := \Hom_{\Mod A}({\cal F}(P), -)$, ${\cal V} :=
\Hom_B(Q,-)$ is inverse to ${\cal U}$, and ${\cal U} \cong {\cal U}'
:= \Hom_A(P,-)$.
\end{proposition}

\begin{proof} Clearly, ${\cal F}(Q) = {\cal F}({\cal U}(A)) =
{\cal L}({\cal F}(A))$. Since ${\cal L}$ is an equivalence,
neglecting grading, we have $A \cong \End_{\Mod A}({\cal F}(A))
\cong \End_{\Mod B}({\cal L}({\cal F}(A)))= \End_{\Mod B}({\cal
F}(Q)) = \End_B(Q)$. This isomorphism maps $a \in A_g$ to ${\cal
L}(r_a)$ $\in \End_B(Q)_g$, where $r_a \in \End_{\Mod A}({\cal
F}(A))$ is the right multiplication by $a$. Thus $A \cong \End_B(Q)$
as graded algebras. Using this isomorphism, we may identify
$\Gr(\End_B(Q))$ with $\Gr A$, and $\Mod(\End_B(Q))$ with $\Mod A$.
So ${\cal M}$ is a functor from $\Mod B$ to $\Mod A$. For any $M \in
\Mod A$, we have $M \cong \Hom_{\Mod A}({\cal F}(A),M) \cong
\Hom_{\Mod B}({\cal L}({\cal F}(A)),{\cal L}(M)) = {\cal M}({\cal
L}(M))$ as $A$-modules. Thus $\eta _M(m):={\cal L}(r_m)$ for all $M
\in \Mod A$ and $m \in M$ defines a natural isomorphism $\eta :
1_{\Mod A} \rightarrow {\cal M}{\cal L} = \Hom_{\Mod B}({\cal
L}({\cal F}(A)),{\cal L}(-)).$ For any $X \in \Gr A$, we have a
graded $A$-module isomorphism $X \cong \Hom_{\Mod A}({\cal F}(A),$
\linebreak ${\cal F}(X)) \cong \Hom_{\Mod B}({\cal L}({\cal
F}(A)),{\cal L}({\cal F}(X))) = \Hom_{\Mod B}({\cal F}({\cal
U}(A)),{\cal F}({\cal U}(X))) = \Hom_B({\cal U}(A),{\cal U}(X)) =
{\cal V}({\cal U}(X))$. Thus $\xi _X(x):={\cal L}(r_x)$ for all $X
\in \Gr A$ and $x \in X_g$ defines a natural isomorphism $\xi :
1_{\Gr A} \rightarrow {\cal V}{\cal U} = \Hom_{\Mod B}({\cal
L}({\cal F}(A)), $ \linebreak ${\cal L}({\cal F}(-))).$ Since
$\xi_{S_g(X)}(x) = {\cal L}(r_x) = \xi_X(x) = S_g(\xi_X)(x)$ for all
$g, h \in G, X \in \Gr A$ and $x \in X_h$, $\xi$ is a graded natural
isomorphism.

Since ${\cal L}$ is an equivalence, ${\cal F}(Q)$ is a finitely
generated projective generator in $\Mod B$. By \cite[Corollary
22.4]{AF}, the functor ${\cal M} = \Hom_{\Mod B}({\cal F}(Q), -)$ is
an equivalence. By Lemma~\ref{lgradehom}, ${\cal M}$ is a graded
equivalence and ${\cal V}=\Hom_B(Q,-)$ is an associated graded
functor of ${\cal M}$. Since ${\cal V}(B) = P$, replacing $A, B,
{\cal L}, {\cal U}, Q, {\cal V}, {\cal M}$ with $B, A, {\cal M},
{\cal V}, P, {\cal U}', {\cal L}'$ respectively, by the same
argument as above, we obtain $B \cong \End_A(P)$ as graded algebras
and ${\cal U}'{\cal V} \cong 1_{\Gr B}$. Therefore ${\cal U}' \cong
{\cal U}'({\cal V}{\cal U}) = ({\cal U}'{\cal V}){\cal U} \cong
{\cal U}$ and $1_{\Gr B} \cong {\cal U}{\cal V}$. So ${\cal U}$ is
inverse to ${\cal V}$.

Since ${\cal M}$ is an equivalence, ${\cal F}(P)$ is a finitely
generated projective generator in $\Mod A$. Thus ${\cal
L}'=\Hom_{\Mod A}({\cal F}(P), -)$ is an equivalence. Furthermore,
${\cal L}'$ is a graded equivalence and ${\cal U}' = \Hom_A(P,-)$ is
an associated graded functor. It follows from classical Morita
theory that $1_{\Mod A} \cong {\cal M}{\cal L}'$ and $1_{\Mod B}
\cong {\cal L}'{\cal M}$. Therefore ${\cal L}' \cong {\cal L}'({\cal
M}{\cal L}) = ({\cal L}'{\cal M}){\cal L} \cong {\cal L}$ and
$1_{\Mod B} \cong {\cal L}{\cal M}$. So ${\cal M}$ is inverse to
${\cal L}$.
\end{proof}

\begin{corollary} Let $A$ and $B$ be graded algebras. Then any
associated graded functor of a graded equivalence $\Mod A
\rightarrow \Mod B$ is a graded equivalence. \end{corollary}

\begin{proposition} \label{passfunctor} Let $A$ and $B$ be graded algebras
and ${\cal U}: \Gr A \rightarrow \Gr B$ a graded equivalence. Then
${\cal U}$ is isomorphic to an associated graded functor of some
graded equivalence $\Mod A \rightarrow \Mod B$. \end{proposition}

\begin{proof} Let ${\cal V}$ be an inverse of ${\cal U}$ and
$P={\cal V}(B)$. Then $P$ is a finitely generated gr-projective
$A$-module. Thus ${\cal F}(P)$ is a finitely generated projective
$A$-module. Since $\oplus_{g \in G}S_g(B)$ is a gr-projective
generator in $\Gr B$ and ${\cal V}$ is a graded equivalent functor,
$\oplus_{g \in G}S_g(P)$ is a gr-projective generator in $\Gr A$.
Hence $A$ is a direct summand of a direct sum of copies of this
generator. It follows that ${\cal F}(P)$ is a finitely generated
projective generator in $\Mod A$. Thus the functor $\Hom_{\Mod
A}({\cal F}(P), -): \Mod A \rightarrow \Mod B$ is an equivalence.

Since ${\cal V}$ is a graded equivalent functor, for each $g \in G$,
there is a unique isomorphism $\phi_g$ such that the following
diagram is commutative:
$$\begin{array}{ccc} \Hom_{\Gr B}(B, S_g(B)) &
\stackrel{{\cal V}}{\longrightarrow} & \Hom_{\Gr A}(P, S_g(P))\\
{\cal F} \downarrow & & \downarrow {\cal F}\\
\Hom_B(B,B)_g & \stackrel{\phi_g}{\longrightarrow} & \Hom_A(P,P)_g.
\end{array}$$ The map $\phi := \oplus_{g \in G} \phi_g$ defines an isomorphism
of graded algebras $\End_B(B) \rightarrow \End_A(P)$. Using this
isomorphism, $\Hom_A(P,X)$ is a graded $B$-module for all $X \in \Gr
A$. By Lemma~\ref{lgradehom}, $\Hom_A(P,-)$ is an associated graded
functor of $\Hom_{\Mod A}({\cal F}(P),-)$. Now it suffices to show
that ${\cal U} \cong \Hom_A(P,-)$. Since ${\cal U} \cong
\Hom_B(B,{\cal U}(-))$, it suffices to show $\Hom_B(B,{\cal U}(-))$
$\cong \Hom_A(P,-)$.

Let $\eta: {\cal V}{\cal U} \rightarrow 1_{\Gr A}$ be a graded
natural isomorphism. Since ${\cal U}$ and ${\cal V}$ are graded
equivalences, for $X \in \Gr A$ and $g \in G$, there is a unique
isomorphism $\xi_{X,g}$ such that the following diagram is
commutative:
$$\begin{array}{ccc} \Hom_{\Gr B}(B,S_g({\cal U}(X))) &
\stackrel{\Hom_{\Gr A}(P, \eta_{S_g(X)})({\cal
V}(-))}{\longrightarrow} &
\Hom_{\Gr A}(P,S_g(X))\\ {\cal F} \downarrow & & \downarrow {\cal F}\\
\Hom_B(B,{\cal U}(X))_g & \stackrel{\xi_{X,g}}{\longrightarrow} &
\Hom_A(P,X)_g.
\end{array}$$
The map $\xi_X := \oplus_{g \in G}\xi_{X,g}: \Hom_B(B,{\cal U}(X))
\rightarrow \Hom_A(P,X)$ is a graded vector space isomorphism. Now
we show that $\xi_X$ is a graded $B$-module morphism, i.e.,
$\xi_{X,hg}(b\psi)=b\xi_{X,g}(\psi)$ for $b \in B_h$ and $\psi \in
\Hom_B(B,{\cal U}(X))_g$. Let $\psi' \in \Hom_{\Gr B}(B,S_g({\cal
U}(X)))$ and $r_b' \in \Hom_{\Gr B}(B,S_h(B))$ such that ${\cal
F}(\psi')=\psi$ and ${\cal F}(r_b')=r_b$. Then $\xi_{X,g}(\psi) =
{\cal F}({\cal V}(\psi')){\cal F}({\eta_{S_g(X)}})$ and
$\xi_{X,hg}(b\psi) = \xi_{X,hg}(r_b \psi) = {\cal F}({\cal V}(r_b'
\psi ')){\cal F} (\eta_{S_{hg}(X)}).$ Thus, for any $p \in P$, we
have \linebreak $(b\xi_{X,g}(\psi))(p) = \xi_{X,g}(\psi)(pb) =
\xi_{X,g}(\psi)(\phi_h(r_b)(p)) = \xi_{X,g}(\psi)({\cal F}({\cal
V}(r_b'))(p)) = ({\cal F}({\cal V}(r_b')){\cal F}({\cal
V}(\psi')){\cal F}({\eta_{S_g(X)}}))(p) = ({\cal F}({\cal V}(r_b'
\psi')){\cal F}({\eta_{S_{hg}(X)}}))(p) =\xi_{X,hg}(b\psi)(p)$,
where the fifth equality holds since $\eta$ is a graded natural
isomorphism. Therefore, $\xi$ is a natural isomorphism. Moreover,
for all $g,h \in G, X \in \Gr A$ and $y \in {\cal U}(S_g(X))_h$,
i.e., $y \in {\cal U}(X)_{hg}$, we have $\xi_{S_g(X)}(r_y) =
\xi_{S_g(X), h}(r_y) = {\cal V}(r_y) \eta_{S_h(S_g(X))} = {\cal
V}(r_y)\eta_{S_{hg}(X)} = \xi_{X, hg}(r_y) = \xi_X (r_y) =
S_g(\xi_X)(r_y)$. Hence $\xi$ is a graded natural isomorphism.
\end{proof}

\begin{theorem} \label{tmoritaequ} Let $A$ and $B$ be graded algebras. Then the
following statements are equivalent:

(1) $A$ and $B$ are graded equivalent;

(2) There is a graded equivalence $\Mod A \rightarrow \Mod B$.

(3) There exists $P \in \Gr A$ such that ${\cal F}(P)$ is a finitely
generated projective generator in $\Mod A$ and $\End_A(P)$ is
isomorphic to $B$ as graded algebras.
\end{theorem}

\begin{proof} $(1) \Leftrightarrow (2)$: It follows from
Proposition~\ref{pmoritaequ} and Proposition~\ref{passfunctor}.

$(2)\Rightarrow (3)$: It follows from the proof of
Proposition~\ref{pmoritaequ}.

$(3)\Rightarrow (2)$. Assume that there exists an object $P$ of $\Gr
A$ such that ${\cal F}(P)$ is a finitely generated projective
generator in $\Mod A$ and $\End_A(P)$ is isomorphic to $B$ as graded
algebras. Then we may identify $\Gr B$ with $\Gr \End_A(P)$, and
$\Mod B$ with $\Mod \End_A(P)$. By Lemma~\ref{lgradehom},
$\Hom_{\Mod A}({\cal F}(P), -)$ is not only an equivalence but also
a graded equivalence.\end{proof}

\subsection{Gr-basic algebras}

\begin{lemma} \label{lendradical} Let $A$ be a finite-dimensional
graded algebra and $M, N \in \gr A$ two gr-indecomposable
$A$-modules satisfying $M \ncong S_g(N)$ for all $g \in G$. Then for
all $\phi \in \Hom_A(M, N)$ and $\psi \in \Hom_A(N, M)$, we have
$\phi\psi \in J_G(\End_{A}(M))$. \end{lemma}

\begin{proof} Let $\phi=\oplus_{g \in G}\phi_g \in \Hom_A(M, N)$ and
$\psi = \oplus_{h \in G}\psi_h \in \Hom_A(N, M)$ with $\phi_g \in
\Hom_A(M, N)_g$ and $\psi_h \in \Hom_A(N, M)_h$. For $g,h \in G$, by
Proposition~\ref{pindecomposable}, $\phi_g \psi_h \in \End_{A}(M)$
is either invertible or nilpotent. If it is invertible then $M \cong
S_{gh}(N)$. It is a contradiction. So it is nilpotent. By the proof
of Proposition~\ref{pgradedlocal}, we have $\phi\psi \in
J_G(\End_A(M))$.
\end{proof}

\begin{lemma} \label{ltopendproj} Let $A$ be a finite-dimensional
graded algebra and $P$ a finite-dimensional gr-indecomposable
gr-projective module. Then as graded algebras
$\End_A(P)/J_G(\End_A(P)) \cong
\End_{A/J_G(A)}(P/J_G(P))$.\end{lemma}

\begin{proof} For $f_g \in
\Hom_{\gr A}(P, S_g(P))$, by \cite[Proposition 2.9.1 iii)]{NV}, we
can define a natural graded morphism $\bar{f}_g \in \Hom_{\gr
A}(P/J_G(P), S_g(P/J_G(P))$. Consequently, we obtain a graded
algebra epimorphism $\phi: \End_A(P) \rightarrow \End_A(P/J_G(P))$,
$\sum_{g \in G} {\cal F}(f_g) \mapsto \sum_{g \in G} {\cal
F}(\bar{f}_g)$.

By Proposition~\ref{pgradedsimple}, we know $P/J_G(P)$ is gr-simple.
It follows from Proposition~\ref{pindecomposable} that
$\End_{A/J_G(A)}(P/J_G(P)) \cong \End_A(P/J_G(P))$ is gr-division.
By Proposition~\ref{pgradedlocal}, we have $\Ker \phi =
J_G(\End_A(P))$. Thus, as graded algebras, $\End_A(P)/J_G(\End_A(P))
\cong \End_{A/J_G(A)}(P/J_G(P))$.
\end{proof}

\begin{proposition} \label{pbasicdef}
Let $A$ be a finite-dimensional graded algebra. Then the following
conditions are equivalent:

(1) Any decomposition of $A$ into gr-indecomposable gr-projective
modules $A = \oplus^r_{i=1}P_i$ satisfies $P_i \ncong
\sigma(g)(P_j)$ for all $g \in G$ and $i \neq j$.

(2) $A/J_G(A)$ is a direct sum of gr-division algebras.
\end{proposition}

\begin{proof} $(1) \Rightarrow (2):$  By
Proposition~\ref{pgradedsimple} and
Proposition~\ref{pindecomposable}, it is enough to show that
$A/J_G(A) \cong \oplus_{i=1}^n\End_{A/J_G(A)}(P_i/J_G(P_i))$ as
graded algebras.

Let $E :=(E_{ij})_{i, j=1}^n$ where $E_{ij} := \Hom_A(P_j, P_i)$ for
$i, j=1, \cdots, n$. Then $E$ is a finite-dimensional $G$-graded
algebra with $E_g := (\Hom_A(P_j, P_i)_g)_{i, j=1}^n$ for all $g \in
G$. Furthermore, $A \cong E$ as $G$-graded algebras.

By Lemma~\ref{lendradical} and Lemma~\ref{ltopendproj}, we can show
that
$$J_G(E) = \left(\begin{array}{cccc}
J_G(E_{11})&E_{12}&\cdots&E_{1n}\\
E_{21}&J_G(E_{22})&\cdots&E_{2n}\\
\cdots&\cdots&\cdots&\cdots\\
E_{n1}&E_{n2}&\cdots&J_G(E_{nn})
\end{array}\right)$$
and as graded algebras$A/J_G(A) \cong \End_A(A)/ J_G(\End_A(A))
\cong E/J_G(E) \cong \oplus_{i=1}^n \End_A(P_i)/J_G(\End_A(P_i))
\cong \oplus_{i=1}^n \End_{A/J_G(A)}(P_i/J_G(P_i))$.

$(2) \Rightarrow (1):$ By Theorem~\ref{tkrullschmidt}, up to
permutation and graded isomorphism, $A$ can be decomposed uniquely
into $A = \oplus_{i=1}^r\oplus_{j=1}^{s_i}S_{g_{ij}}(P_i)$, where
$P_i$'s are gr-indecomposable gr-projective $A$-modules such that
$P_i \ncong S_g(P_j)$ for all $g \in G$ and $i\neq j$. Using the
same argument as above, we can obtain $A/J_G(A) \cong
\oplus_{i=1}^r\End_A(\oplus_{j=1}^{s_i}S_{g_{ij}}(P_i))/
J_G(\End_A(\oplus_{j=1}^{s_i}S_{g_{ij}}(P_i))).$ If $s_i>1$ for some
$i$ then it is easy to find a nonzero nilpotent homogeneous element
in $A/J_G(A)$. However, this is impossible in the direct sum of
gr-division algebras. \end{proof}

\begin{definition} We say a finite-dimensional graded algebra $A$ is
\mbox{\bf gr-basic} if it satisfies the equivalent conditions of
Proposition~\ref{pbasicdef}. \end{definition}

Let $A$ be a finite-dimensional gr-basic graded algebra and $\{P_i |
1 \leq i \leq n\}$ a complete set of gr-indecomposable gr-projective
$A$-modules up to graded isomorphism and shift. A graded algebra $B$
is called a {\it congruence} of $A$ if it is isomorphic to
$\End_A(\oplus_{i=1}^nS_{g_i}(P_i))$ for some $(g_1,...,g_n) \in
G^n$.

\begin{theorem} \label{treducedtobasic} Any finite-dimensional graded algebra $A$ is
graded equivalent to a gr-basic algebra $B$, which is unique up to
congruence.
\end{theorem}

\begin{proof} By Theorem~\ref{tkrullschmidt},
up to permutation and graded isomorphism, $A$ can be decomposed
uniquely into $A = \oplus_{i=1}^n\oplus_{j=1}^{r_i}S_{g_{ij}}(P_i)$,
where $P_i$'s are gr-indecomposable gr-projective $A$-modules such
that $P_i \ncong S_g(P_j)$ for all $g \in G$ and $i \neq j$. Let $P
= \oplus_{i=1}^nS_{g_{i1}}(P_i)$. Note that ${\cal F}(A)=
\oplus_{i=1}^n\oplus_{j=1}^{r_i}{\cal F}(S_{g_{ij}}(P_i))$ is a
finitely generated projective generator in $\Mod A$, so is ${\cal
F}(P)$. It follows from Theorem~\ref{tmoritaequ} that $A$ is graded
equivalent to $B : = \End_A(P)$. Since $B/J_G(B) \cong
\oplus_{i=1}^n\End_{A/J_G(A)}(S_{g_{i1}}(P_i)/J_G(S_{g_{i1}}(P_i)))$
as graded algebras, by Proposition~\ref{pbasicdef}, $B$ is gr-basic.

Assume that $A$ is also graded equivalent to another gr-basic graded
algebra $C$. Then $B$ is graded equivalent to $C$. Suppose that
${\cal H}: \Gr B \rightarrow \Gr C$ is a graded equivalence, and
$B=\oplus_{i=1}^mQ_i$ and $C=\oplus_{j=1}^nR_j$ are the
decompositions of $B$ and $C$ into gr-indecomposable gr-projective
modules respectively. Then $m=n$. Moreover, by
Proposition~\ref{pgradedprojective}, $B \cong
\End_B(\oplus_{i=1}^nQ_i) \cong \End_C(\oplus_{i=1}^n{\cal H}(Q_i))
\cong \End_C(\oplus_{i=1}^nS_{g_i}(R_i))$ for some $(g_1,...,g_n)
\in G^n$ , which is a congruence of $\End_C(\oplus_{i=1}^nR_i) \cong
C$. \end{proof}

\section{Constructions of superalgebras}

In this section, we shall consider the nice constructions of
finite-dimensional superalgebras. Let $Q$ be a quiver with vertex
set $Q_0$ and arrow set $Q_1$. Then any map $\mathrm{deg}: Q_1
\rightarrow \mathbb{Z}_2$ induces a $\mathbb{Z}_2$-grading on the
path algebra $KQ$. Thus $KQ/I$ is a finite-dimensional gr-basic
superalgebra for each admissible graded ideal $I$ of $KQ$. It is
natural to ask whether all finite-dimensional gr-basic superalgebras
can be obtained in this way? The answer is NO! Indeed, exactly all
elementary superalgebras can be obtained in this way. For the
knowledge of quivers and their representation theory, we refer to
\cite{ARS}.

\subsection{Elementary superalgebras}

\indent\indent An {\it elementary superalgebra} $A$ is a
finite-dimensional superalgebra such that $A/J_G(A) \cong K \times
\cdots \times K$. Clearly, an elementary superalgebra is gr-basic.

A {\it superquiver} is an oriented diagram whose vertices are either
white vertices $\circ$ or black vertices $\bullet$, and whose arrows
are either solid arrows $\unitlength=1mm
\begin{picture}(8,1) \put(0,1){\vector(1,0){8}}
\end{picture}$ or dotted arrows
$\unitlength=1mm \begin{picture}(8,1)
\multiput(0,1)(1,0){8}{\circle*{0.1}} \put(8,1){\vector(1,0){0}}
\end{picture}$.
We say a superquiver is \textit{1-color} if all its vertices are of
the same color, and \textit{2-color} otherwise. An {\it elementary
superquiver} is just a superquiver with only white vertices. The
{\it underlying quiver} $\underline{Q}$ of a superquiver $Q$ is the
quiver obtained from $Q$ by changing all black vertices into white
ones and all dotted arrows into solid ones respectively. The {\it
underlying diagram} $\overline{Q}$ of a quiver or a superquiver $Q$
is the diagram obtained from $Q$ by changing all black vertices into
white ones and all arrows into edges.

For an elementary superquiver $Q$, we define its {\it path
superalgebra} $KQ$ to be the superalgebra which has all paths as a
$K$-basis, whose multiplication is given by the concatenation of
paths, and the $\mathbb{Z}_2$-grading on $KQ$ is given by
$\mathrm{deg}(\circ)=\mathrm{deg}(\unitlength=1mm
\begin{picture}(8,1) \put(0,1){\vector(1,0){8}}
\end{picture}): = \bar{0}$ and $\mathrm{deg}(\unitlength=1mm \begin{picture}(8,1)
\multiput(0,1)(1,0){8}{\circle*{0.1}} \put(8,1){\vector(1,0){0}}
\end{picture}): = \bar{1}$.

A set $\{e_1, \cdots, e_n\}$ of idempotents in a $G$-graded algebra
$A$ is said to be \textit{orthogonal} if $e_ie_j=0$ for all $i\neq
j$. A nonzero idempotent $e$ is said to be \textit{primitive} if $e$
cannot be written as a sum of two nonzero orthogonal idempotents. We
say a set $\{e_1, \cdots, e_n\}$ of orthogonal idempotents in $A$ is
\textit{complete} if $1=e_1+\cdots +e_n$.

\begin{proposition}\label{pidempotent}
Each finite-dimensional graded algebra $A$ has a complete set of
degree $e$ primitive orthogonal idempotents.
\end{proposition}

\begin{proof}
Let $A=P_1 \oplus \cdots \oplus P_n$ be a decomposition of $A$ into
a direct sum of gr-indecomposable gr-projective modules. Then there
exist $e_i \in (P_i)_e$ such that $1=e_1+\cdots+e_n$. Thus $\{e_1,
\cdots, e_n\}$ is a complete set of primitive orthogonal idempotents
of $A$: Indeed, for any $x_i \in P_i$, we have
$x_i=x_ie_1+\cdots+x_ie_i+\cdots+x_ie_n$. Thus $x_ie_j=0$ for $j\neq
i$ and $x_i=x_ie_i$. In particular, we have $e_i \neq 0, e_i^2=e_i$
and $e_ie_j=0$ for $i\neq j$. Moreover, $e_i$'s are primitive for
all $i$: Assume on the contrary that $e_i=e'+e''$ with $e', e''$
nonzero orthogonal idempotents. Then $P_i=Ae_i=Ae' \oplus Ae''$. It
is a contradiction. \end{proof}

\begin{theorem} An elementary superalgebra $A$ is graded
isomorphic to $KQ/I$ where $KQ$ is the path superalgebra of an
elementary superquiver $Q$ and $I$ is an admissible graded ideal of
$KQ$. If further $A$ is gr-hereditary then $A \cong KQ$ for an
acyclic elementary superquiver $Q$. \end{theorem}

\begin{proof} The first part of the theorem can be proved analogous
to \cite[Theorem III.1.9]{ARS} by using
Proposition~\ref{pidempotent}. The second part of the theorem can be
proved analogous to \cite[Proposition III.1.13]{ARS}.
\end{proof}

\subsection{Graded species}

\indent\indent We have already provided the nice and easy
constructions above for elementary superalgebras. To deal with
general cases, we introduce graded species and their
representations, which are interesting in their own rights.

\begin{definition} A \mbox{\bf $G$-graded $K$-species} or \mbox{\bf graded species}
$\mathscr{S} =$ \linebreak $(D_i, \, _jM_i)_{i, j \in I}$ with $I :=
\{ 1, 2, \ldots, n \}$ is a collection of finite-dimensional
gr-division algebras $D_i$ and finite-dimensional graded
$D_j$-$D_i$-bimodules $_jM_i$ such that $K$ operates on $_jM_i$
centrally, i.e., $km=mk$ for all $k \in K$ and $m \in \, _jM_i$. A
$\mathbb{Z}_2$-graded $K$-species is called a {\bf superspecies}.
\end{definition}


Let $A$ be a finite-dimensional gr-basic graded algebra. Then
$A/J_G(A) = \oplus_{i=1}^n D_i$ where $D_i$'s are gr-division
algebras, and $J_G(A)/J_G(A)^2 = \oplus_{i,j \in I}{_jM_i}$ where
$_jM_i$ is a graded $D_j$-$D_i$-bimodule for all $i,j \in I$. We
call $\mathscr{S}_A = (D_i, \, _jM_i)_{i, j \in I}$ the {\it graded
species of} $A$.

A {\it representation} $V=(V_i, \, _j\phi_i)_{i, j \in I}$ of a
graded species $\mathscr{S}=$ \linebreak $(D_i, \, _jM_i)_{i, j \in
I}$ is a collection of graded $D_i$-modules $V_i$ and graded
$D_j$-module morphisms $_j\phi_i: \, _jM_i \otimes_{D_i} V_i
\rightarrow V_j$ for all $i, j \in I$. We say a representation
$V=(V_i, \, _j\phi_i)_{i, j \in I}$ is {\it finite-dimensional} if
all $V_i$ are finite-dimensional. Let $V=(V_i, \, _j\phi_i)_{i, j
\in I}$ and $W=(W_i, \, _j\psi_i)_{i, j \in I}$ be representations
of $\mathscr{S}$. Then the {\it direct sum} of the representations
$V$ and $W$ is the representation $V \oplus W := (V_i \oplus W_i,
\mbox{diag}\{\, _j\phi_i, \,_j\psi_i\})_{i, j \in I}$. We say a
representation of a graded species is {\it indecomposable} if it is
not the direct sum of two nonzero representations. A {\it morphism}
$\az = (\az_i)_{i \in I} : V \rightarrow W$ is a collection of
graded $D_i$-module morphisms $\az_i: V_i \rightarrow W_j$ such that
$(1 \otimes \az_i) \,_j\psi_i = \, _j\phi_i \,\, \az_j$ for all $i,
j \in I$. All representations of $\mathscr{S}$ and all morphisms
between them form an abelian category $\Rep \mathscr{S}$. We denote
by $\rep \s$ the full subcategory of $\Rep \mathscr{S}$ consisting
of all finite-dimensional representations.

Let $R$ be a graded algebra and $M$ a graded $R$-$R$-bimodule. With
the pair $(R, M)$ we associate a {\it graded tensor algebra} $T(R,
M) = \oplus_{i=0}^{\infty}M^{\otimes_R i}$, where $M^{\otimes_R
0}=R$, $M^{\otimes_R 1}= M$, $M^{\otimes_R i}$ is the $i$-folds
tensor product $M \otimes_R \cdots \otimes_R M$, and the
multiplication is induced by the tensor product. We say a graded
ideal $\i$ of $T(R, M)$ is {\it admissible} if $(\oplus_{i \geq
1}M^{\otimes_R i})^t \subseteq \i \subseteq (\oplus_{i \geq
1}M^{\otimes_R i})^2$ for some $t \ge 2$.

With a graded species $\mathscr{S}= (D_i, \, _jM_i)_{i, j \in I}$ we
associate a {\it graded tensor algebra} $T(\mathscr{S}) = T(R,M)$,
where $R  = \oplus_{i \in I}D_i$ and $M = \oplus_{i, j \in I} \,
_jM_i$. A graded species $\mathscr{S}$ is said to be {\it acyclic}
if its graded tensor algebra $T(\mathscr{S})$ is finite-dimensional.

\begin{proposition} \label{pspeciestensor}
Let $\mathscr{S}= (D_i, \, _jM_i)_{i, j \in I}$ be a graded species.
Then $\Rep \mathscr{S}$ and $\rep \mathscr{S}$ are equivalent to
$\Gr T(\mathscr{S})$ and $\gr T(\mathscr{S})$ respectively.
\end{proposition}

\begin{proof} It suffices to add ``graded'' or ``gr-'' at appropriate
places of the proof of \cite[Proposition 10.1]{DR}. \end{proof}

\begin{remark} If $\s$ is an acyclic graded species then by
Theorem~\ref{tkrullschmidt} and Proposition~\ref{pspeciestensor},
$\rep \mathscr{S}$ is a Krull-Schmidt category (ref.
\cite[p.52]{R1099}). \end{remark}

\subsection{Construction of finite-dimensional superalgebras}

\indent\indent Let $A$ be a finite-dimensional superalgebra. Then
both its opposite algebra $A^{\op}$ and its enveloping algebra
$A^e=A \ot_K A^{\op}$ are finite-dimensional superalgebras with
natural grading.

\begin{proposition}\label{pseparable}
Let $A$ be a finite-dimensional superalgebra and $\phi : A^e
\rightarrow A$ the graded $A^e$-module epimorphism given by $\phi(x
\ot y^o)=xy$. Then the following conditions are equivalent:

(1) $A$ is a gr-projective $A^e$-module;

(2) There exists a degree $\bar{0}$ element $\epsilon$ of the graded
$A$-$A$-bimodule $A \ot_K A$ (isomorphic to $A^e$) such that
$\phi(\epsilon)=1$ and $a \epsilon=\epsilon a$ for all $a\in A$.
\end{proposition}

\begin{proof} (1) $\Rightarrow$ (2): Since $A$ is a gr-projective
$A^e$-module, there is a graded $A^e$-module morphism $\psi : A
\rightarrow A^e$ such that $\psi \phi = \mbox{id}_A$. Thus $\epsilon
:= \psi(1_A)$ is as required.

(2) $\Rightarrow$ (1): The map $\psi : A \rightarrow A^e, a \mapsto
a \epsilon = \epsilon a$, is a graded $A^e$-module morphism and
$\psi \phi = \mbox{id}_A$.
\end{proof}

A finite-dimensional superalgebra is said to be \textit{super
separable} if it satisfies the equivalent conditions in
Proposition~\ref{pseparable}.

\begin{lemma} \label{ldirectsumsuperseparable} The direct product of
finitely many super separable algebras is super separable.
\end{lemma}

\begin{proof} If $A_1$ and $A_2$ are super separable then there are
degree $\bar{0}$ elements $\epsilon_i \in A_i \otimes _K A_i$ such
that $\phi_i(\epsilon_i)=1_{A_i}$ and $a_i \epsilon_i = \epsilon_i
a_i$ for all $a_i \in A_i$ and $i=1,2$, where $\phi_i : A_i^e
\rightarrow A_i, a_i \otimes b_i^o \mapsto a_ib_i$. Note that there
is a natural graded $(A_1 \times A_2)^e$-module morphism $\varphi :
(A_1 \otimes _K A_1) \times (A_2 \otimes _K A_2) \rightarrow (A_1
\times A_2) \otimes _K (A_1 \times A_2), (a_1 \otimes b_1, a_2
\otimes b_2) \mapsto (a_1,0) \otimes (0, a_2) + (0, b_1) \otimes (0,
b_2)$. The element $\epsilon := \varphi(\epsilon_1,\epsilon_2)$ is
as required. \end{proof}

\begin{lemma} \label{lgradeddivision} Let $A = \oplus_{g \in G}
A_g$ be a gr-division algebra. Then $A_g = A_e \varepsilon$ for some
$\varepsilon \in A_g$.
\end{lemma}

\begin{proof} If $A_g = 0$ then we may take $\varepsilon=0$. If $A_g \neq
0$ then we may take any $\varepsilon \in A_g \backslash \{0\}$.
Since $\varepsilon$ is invertible, we have $\varepsilon^{-1} \in
A_{g^{-1}}$. Thus $x\varepsilon^{-1} \in A_e$ for all $x \in A_g$.
Hence $x \in A_e \varepsilon$ and $A_g = A_e \varepsilon$.
\end{proof}

\begin{proposition} \label{psuperdivision} Let $K$ be an algebraically
closed field and $A$ a finite-dimensional gr-division superalgebra.
Then $A$ is either $K$ with trivial grading or $D$ with $D_{\bar{0}}
= K$, $D_{\bar{1}} = K \varepsilon$ and $\varepsilon^2=1$.
\end{proposition}

\begin{proof} By Lemma~\ref{lgradeddivision}, we have $A = A_{\bar{0}}
\oplus A_{\bar{0}} \varepsilon$ for some $\varepsilon \in
A_{\bar{1}}$. Since $A_{\bar{0}}$ is a finite-dimensional division
algebra and $K$ is algebraically closed, we have $A_{\bar{0}} = K$
by \cite[Lemma 3.5]{P}. Thus $A = K \oplus K \varepsilon$ with
$\varepsilon^2 \in K$. If $\varepsilon = 0$ then $A=K$. If
$\varepsilon \neq 0$ then we may assume that $\varepsilon^2=1$,
since $K$ is algebraically closed. \end{proof}

\begin{proposition}\label{pquasifree}
Let $K$ be an algebraically closed field of $\mathrm{char} K \neq 2$
and $A$ a finite-dimensional gr-basic gr-semisimple superalgebra.
Then $A$ is super separable.
\end{proposition}

\begin{proof} By assumption and Proposition~\ref{psuperdivision}, $A$ is isomorphic to $\oplus_{i=1}^nD_i$
where $D_i=K$ or $D$ for all $i$. Clearly, $K$ is super separable.
In order to show that $D$ is super separable, it suffices to take
$\epsilon=\frac{1}{2}(1 \otimes 1+\varepsilon \ot \varepsilon)$
which satisfies the second condition in
Proposition~\ref{pseparable}. It follows from
Lemma~\ref{ldirectsumsuperseparable} that $A$ is super separable.
\end{proof}

A \textit{square zero super extension} of a superalgebra $R$ by a
graded $R$-$R$-bimodule $M$ is a superalgebra $E$, together with a
superalgebra epimorphism $\phi: E \rightarrow R$ such that $\Ker
\phi$ is a square zero graded ideal, and a graded $R$-$R$-bimodule
isomorphism of $M$ with $\Ker \phi$ (cf. \cite[\S 9.3]{W}). We say
two square zero super extensions $E$ and $E'$ of $R$ by $M$ are
\textit{equivalent} if there is a superalgebra isomorphism $\psi: E
\rightarrow E'$ such that the following diagram commutes:
$$\begin{array}{ccccc}
  M&\rightarrowtail&E&\twoheadrightarrow&R \\
  \parallel&& \psi\!\!\downarrow&& \parallel\\
  M&\rightarrowtail&E'&\twoheadrightarrow&R \\
\end{array}$$

\begin{proposition} \label{squarezeroextension} Let $R$ be a superalgebra and $M$
a graded $R$-$R$-bimodule. Then the equivalence classes of square
zero super extensions of $R$ by $M$ are in 1-1 correspondence with
the elements of $\Ext_{\Gr R^e}^2(R, M)$.
\end{proposition}

\begin{proof} Analogous to the proof of \cite[Theorem
9.3.1]{W}. It suffices to note that all boundary maps in the bar
resolution of $R$ are graded (ref. \cite{Lo}). \end{proof}

We say that a superalgebra $R$ is \textit{super quasifree} if for
every square zero super extension $M \rightarrowtail E
\stackrel{\phi} \twoheadrightarrow T$ of a superalgebra $T$ by a
graded $T$-$T$-bimodule $M$ and every superalgebra homomorphism
$\psi: R \rightarrow T$, there exists a superalgebra homomorphism
$\varphi: R \rightarrow E$ lifting $\psi$, i.e., $\varphi\phi =
\psi$.

\begin{proposition} \label{pquasi-free}
A superalgebra $R$ is super quasifree if and only if \linebreak
$\Ext_{\Gr R^e}^2(R, M)=0$ for all graded $R$-$R$-bimodule $M$.
\end{proposition}

\begin{proof} Owing to Proposition~\ref{squarezeroextension}, we
need do little modification on the proof of \cite[Proposition
9.3.3]{W}. \end{proof}

\begin{theorem} \label{tWedderburn} {\bf (Super version of Wedderburn's principal
theorem)} Let $K$ be an algebraically closed field of $\ch K\neq 2$
and $A$ a finite-dimensional gr-basic superalgebra. Then there is a
sub-superalgebra $S$ of $A$ such that $A=S \oplus J_{\z_2}(A)$.
\end{theorem}

\begin{proof}
Since $A$ is gr-basic, we know $R := A/J_{\z_2}(A) =
\oplus_{i=1}^nD_i$ is a gr-basic gr-semisimple superalgebra. By
Proposition~\ref{pquasifree}, $R$ is super separable. Thus
$\Ext_{\Gr R^e}^n(R,M)=0$ for all $n \geq 1$ and $M \in \Gr R^e$. It
follows from Proposition~\ref{pquasi-free} that $R$ is super
quasifree. Thus the superalgebra identity $R \rightarrow R$ can be
lifted successively to $R \rightarrow A/J_{\z_2}^2(A)$, $R
\rightarrow A/J_{\z_2}(A)^3, \cdots$. Since $J_{\z_2}(A)^m=0$ for
some $m$, it is eventually lifted to $\phi : R \rightarrow
A/J_{\z_2}(A)^m=A$. Hence there is a sub-superalgebra $S=\mathrm{Im}
\phi$ of $A$ such that $A=S \oplus J_{\z_2}(A)$.
\end{proof}

\begin{lemma}\label{lextension} Let $R$ be a graded
algebra, $M$ a graded $R$-$R$-bimodule, and $A$ a graded algebra.
Suppose $\varphi: R \oplus M \rightarrow A$ is a map satisfying:

(1) $\varphi|_R: R \rightarrow A$ is a graded algebra homomorphism;

(2) $\varphi|_M: M \rightarrow A$ is a graded $R$-$R$-bimodule
morphism when $A$ is viewed as a graded $R$-$R$-bimodule via
$\varphi|_R: R \rightarrow A$.

Then there is a unique graded algebra homomorphism $\tilde{\varphi}:
T(R, M)\rightarrow A$ such that $\tilde{\varphi}|_{R \oplus
M}=\varphi$.\end{lemma}

\begin{proof} First of all, we have a unique graded $R$-$R$-bimodule morphism $\varphi_i:
M^{\otimes_R i} \rightarrow A$ such that
$\varphi_i(x_1\otimes\cdots\otimes x_i)=\varphi(x_1)\cdots
\varphi(x_i)$ for all $i$ and $x_1,...,x_i \in M$. Thus the graded
algebra homomorphism $\tilde{\varphi}: T(R, M) \rightarrow A$ given
by $\tilde{\varphi}(\sum_{i=0}^{\infty}y_i)=\sum_{i=0
}^{\infty}\varphi_i(y_i)$ with $y_i\in M^{\otimes_R i}$ is as
required, which is uniquely determined by $\varphi$. \end{proof}

The following lemma is a graded analog of \cite[Theorem I]{EN}.

\begin{lemma}\label{lhereditary} Let $A$ be a finite-dimension graded algebra.
If $\i \subseteq J_G(A)^2$ is a graded ideal in $A$ such that
$\mbox{gr-gl.dim}\,(A/\i)\le 1$ then $\i=0$.
\end{lemma}

\begin{proof} Clearly, $\i/\i J_G(A) \rightarrowtail J_G(A)/\i J_G(A) \twoheadrightarrow J_G(A)/\i $
is exact in $\Gr (A/\i)$. It follows from $\i \subseteq J_G(A)^2$
that $\top_G(J_G(A)/\i J_G(A)) \cong \top_G(J_G(A)/\i)$. Since
$\mbox{gr-gl.dim} (A/\i) \le 1$ and $J_G(A/\i)=J_G(A)/\i$, we know
$J_G(A)/\i$ is a gr-projective $A/\i$-module. By
Corollary~\ref{cprojcovertop}, the lift of the graded morphism
$J_G(A)/\i \twoheadrightarrow \top_G(J_G(A)/\i)
\stackrel{\cong}{\rightarrow} \top_G(J_G(A)/\i J_G(A))$ to
$J_G(A)/\i \rightarrow J_G(A)/\i J_G(A)$ is a gr-projective cover of
$J_G(A)/\i J_G(A)$. Comparing $K$-dimensions we obtain $J_G(A)/\i
J_G(A)=J_G(A)/\i$, \linebreak i.e., $\i J_G(A)=\i$. Thus $\i=0$.
\end{proof}

\begin{theorem} \label{tpresentation}
Let $K$ be an algebraically closed field of $\ch K\neq 2$, $A$ a
finite-dimensional gr-basic superalgebra and $\s_A=(D_i, \,
_jM_i)_{i,j \in I}$ its superspecies. Then $A$ is graded isomorphic
to $T(\mathscr{S}_A)/ \i$ where $\i$ is an admissible graded ideal
of $T(\s_A)$. Furthermore, if $A$ is gr-hereditary then $\s_A$ is
acyclic and $A$ is graded isomorphic to $T(\s_A)$.
\end{theorem}

\begin{proof} Due to $\ch K\neq 2$ and \cite[Theorem 4.4]{CM},
we may let $J:=J_{\z _2}(A)=J(A)$, $R:=A/J$, and $M:=J/J^2$. By
definition, $T(\s_A)=T(R, M)$. It follows from
Theorem~\ref{tWedderburn} that $A=S \oplus J$ for some
sub-superalgebra $S$ of $A$. Thus there is a natural superalgebra
isomorphism $\varphi_0 : R=A/J \rightarrow S$.

Forgotten grading, $S \cong A/J$ is isomorphic to the direct product
of some copies of $K$, so $S$ is separable and $S^e$ is semisimple.
Owing to $\ch K\neq 2$, applying \cite[Theorem 2.2]{CM},
\cite[Theorem 1.3 and Theorem 1.4]{RR} or \cite[Theorem
III.4.4]{ARS}, and \cite[Theorem 3.5]{CM} in turn, we obtain
$\mbox{gr-gl.dim} S^e = \mbox{gl.dim} S^e \# K\mathbb{Z}_2^* =
\mbox{gl.dim} S^e \# K\mathbb{Z}_2^* \ast \mathbb{Z}_2 =
\mbox{gl.dim} S^e = 0$, i.e., $S^e$ is a gr-semisimple superalgebra.
Hence the exact sequence of graded $S$-$S$-bimodules $J^2
\rightarrowtail J \twoheadrightarrow J/J^2$ splits. So $J=N \oplus
J^2$ for some graded $S$-$S$-bimodule $N$ of $J$. Thus there is a
natural graded $R$-$R$-bimodule isomorphism $\varphi_1 : M=J/J^2
\rightarrow N$, where the graded $R$-$R$-bimodule structure of $M$
is induced by $\varphi_0$. Denote by $N^i$ the set consisting of all
finite sum of the products of $i$ elements in $N$. Then $N^i$ is a
graded $S$-$S$-bimodule. By induction we obtain $J^i=N^i \oplus
J^{i+1}$ for all $i \geq 1$. Suppose $J^{m+1}=0$ for some $m \geq
0$. Then $A=S \oplus N \oplus N^2 \oplus \cdots \oplus N^m$.

The natural injection $\varphi=\varphi_0 \oplus \varphi_1: R \oplus
M \rightarrow A$ satisfies the requirements of
Lemma~\ref{lextension}, thus there is a superalgebra homomorphism
$\tilde{\varphi}: T(R, M) \rightarrow A$ such that
$\tilde{\varphi}|_{R \oplus M} = \varphi$. Since $A=S \oplus N
\oplus N^2 \oplus \cdots \oplus N^m$ and $\varphi_0 , \varphi_1$ are
graded epimorphisms, $\tilde{\varphi}$ is a graded epimorphism and
$A \cong T(R,M)/\Ker \tilde{\varphi}$. Note that
$\tilde{\varphi}(M^{\otimes_R i})\subseteq J^i \subseteq J^2$ for $i
\ge 2$ and $\tilde{\varphi}|_{R \oplus M}$ is an isomorphism, so
$\Ker \tilde{\varphi} \subseteq \oplus_{i \ge 2}M^{\otimes_R i}$.
Since $J^t=0$ for some $t$, we have $\oplus_{i \ge t}M^{\otimes_R i}
\subseteq \Ker\,\tilde{\varphi}$. Thus $A \cong T(R,M)/ \i$ where
$\i := \Ker \tilde{\varphi}$ is an admissible graded ideal of
$T(\s_A)$.

Now assume that $A$ is a finite-dimensional gr-basic gr-hereditary
superalgebra. By above proof we have $A \cong T(\mathscr{S}_A)/\i$
where $\i$ is an admissible graded ideal of $T(\s_A)$. Thus
$T(\mathscr{S}_A)$ must be finite-dimensional: Indeed, by
Proposition~\ref{pidempotent} and its proof, we may assume that
$1=\sum_{i=1}^ne_i$ is a decomposition of $1$ into degree $\bar{0}$
primitive orthogonal idempotents such that $e_iRe_i=D_i$ and
$e_jMe_i=\, _jM_i$ for all $i,j \in I=\{1,2,...,n\}$. If
$T(\mathscr{S}_A)$ is infinite-dimensional then there is $s \gg 0$
such that $M^{\otimes_R s} \neq 0$. Thus there are $j_1,...,j_s \in
I$ such that $_{j_r}M_{j_{r+1}} \neq 0$ for all $1 \leq r \leq s$
where $j_{s+1}=j_1$. Therefore $e_{j_r}Je_{j_{r+1}} \neq 0$ for all
$1 \leq r \leq s$. Suppose $e_{j_r}xe_{j_{r+1}} \in
e_{j_r}Je_{j_{r+1}} \backslash \{0\}$ for some $x \in J$ and all $1
\leq r \leq s$. Then the right multiplication by
$e_{j_r}xe_{j_{r+1}}$ defines a graded morphism $\lambda_r$ from the
gr-indecomposable gr-projective module $Ae_{j_r}$ to the
gr-indecomposable gr-projective module $Ae_{j_{r+1}}$ for each $r$.
Since $A$ is gr-hereditary, $\lambda_r$ must be a graded
monomorphism. It follows from $x \in J$ that $\lambda_r$ is not
surjective. Thus there is a proper injection chain $Ae_{j_1}
\hookrightarrow Ae_{j_2} \hookrightarrow \cdots \hookrightarrow
Ae_{j_s} \hookrightarrow Ae_{j_1}$. It is a contradiction. Hence
$T(\mathscr{S}_A)$ is finite-dimensional, i.e., $\s_A$ is acyclic.
By Lemma~\ref{lhereditary}, we have $\i =0$ and $T(\s_A)$ is graded
isomorphic to $A$.
\end{proof}

\begin{remark} Now we can construct all finite-dimensional superalgebras
over an algebraically closed field $K$ of $\ch K\neq 2$:

(1) Firstly, by Theorem~\ref{tpresentation}, we can construct all
finite-dimensional gr-basic superalgebras by giving a superspecies
$\mathscr{S}=(D_i,\, _jM_i)_{i,j \in I}$ and an admissible graded
ideal $\i$ of $T(\mathscr{S})$. Note that $D_i$ is either $K$ or $D$
and $\, _jM_i$ is nothing but a graded $D_j \otimes_K
D_i^{\op}$-module. In the case of $D_i = D_j = K$, $\, _jM_i$ is
just a super vector space. In the case of $D_i=K$ and $D_j = D$, $\,
_jM_i$ is just a left gr-free $D$-module. In the case of $D_i=D$ and
$D_j = K$, $\, _jM_i$ is just a right gr-free $D$-module. In the
case of $D_i = D_j = D$, $\, _jM_i$ is just a graded $D^e$-module,
equivalently, a direct sum of two gr-free $D$-modules, since $D^e =
D \otimes_K D^{\op} = D \otimes_K D \cong D \oplus D$. Indeed, the
super vector space map $\phi : D \otimes_K D \rightarrow D \oplus D$
given by $1 \otimes 1 \mapsto (1,1), 1 \otimes \varepsilon \mapsto
(\varepsilon,\varepsilon), \varepsilon \otimes 1 \mapsto
(\varepsilon,-\varepsilon), \varepsilon \otimes \varepsilon \mapsto
(1,-1),$ is a superalgebra isomorphism.

(2) Secondly, by Theorem~\ref{tmoritaequ}, we can construct all
finite-dimensional superalgebras as $\End_A(P)$, where $A$ is a
finite-dimensional gr-basic superalgebra with a decomposition
$A=\oplus_{i=1}^nP_i$ of gr-indecomposable gr-projectives and
$P=\oplus_{i=1}^n\oplus_{j=1}^{r_i}S_{g_{ij}}(P_i)$ for some $r_i
\geq 1$ and $g_{ij} \in \mathbb{Z}_2$.

\end{remark}

\section{Graded representation type}

\indent\indent In order to define the graded representation types of
acyclic superspecies, we shall introduce the graded representation
types of finite-dimensional finite group graded algebras at first.
Then we shall obtain the graded version of Drozd's theorem. In this
section, the underlying field $K$ is assumed to be algebraically
closed.

Let $K \langle x,y \rangle$ (resp. $K[x]$) be the free associative
algebra in two variables $x, y$ (resp. one variable $x$). Denote by
$\mod K \langle x,y \rangle$ the category of finite-dimensional left
$K \langle x,y \rangle$-modules. The definition of representation
types of a finite-dimensional $K$-algebra is well-known (ref.
\cite{CB,Dr}). Now we define the graded representation types of
graded $K$-algebras (compare with \cite{HKN, KP}).

\begin{definition} We say that a finite-dimensional $G$-graded $K$-algebra
$A$ is \mbox{\bf gr-representation-finite} provided it has only
finitely many isomorphism classes of gr-indecomposable $A$-modules.
\end{definition}

\begin{definition} We say that a finite-dimensional $G$-graded $K$-algebra
$A$ is \mbox{\bf gr-tame} if for each dimension $d>0$, there are a
finite number of $A$-$K[x]$-bimodules $M_i$ which are $G$-graded as
left $A$-modules and free as right $K[x]$-modules such that every
gr-indecomposable $A$-module of dimension $d$ is isomorphic to $M_i
\otimes_{K[x]} N$ for some $i$ and some simple $K[x]$-module
$N$.\end{definition}

\begin{definition} We say that a finite-dimensional $G$-graded
$K$-algebra $A$ is \mbox{\bf gr-wild} if there is a finitely
generated $A$-$K \langle x, y \rangle$-bimodule $M$ which is
$G$-graded as a left $A$-module and free as a right $K \langle x,y
\rangle$-module and such that the functor $M\otimes_{K \langle x,y
\rangle} -$ from $\mod K \langle x, y \rangle$ to $\gr A$ preserves
indecomposability and isomorphism classes.\end{definition}

\begin{proposition} \label{preptypecoincide} A finite-dimensional
graded algebra $A$ is gr-representation-finite (resp. gr-tame,
gr-wild) if and only if the smash product $A\# KG^*$ is
representation-finite (resp. tame, wild). \end{proposition}

\begin{proof} Let $M$ be a left $A \# KG^*$-module. Then $M$ is also
a $G$-graded left $A$-module defined by $M_g := (1 \# p_g)M$ and $am
:= (a \# 1)m$ for all $g \in G, a \in A$ and $m \in M$. Conversely,
let $M$ be a $G$-graded left $A$-module. Then $M$ is also a left $A
\# KG^*$-module defined by $(a \# p_g)m=am_g$ for all $a \in A$, $g
\in G$ and $m \in M$ (ref. \cite[Section 2]{CM}). By \cite[Theorem
2.2]{CM}, we have that $M$ is an indecomposable $A \# KG^*$-module
if and only if it is a gr-indecomposable $A$-module, and $M$ and $N$
are isomorphic as $A \# KG^*$-modules if and only if they are graded
isomorphic as graded $A$-modules. Furthermore, the
$A$-$K[X]$-bimodules $M_i$ and $A$-$K \langle x, y \rangle$-bimodule
$M$ in the definitions of gr-tameness and gr-wildness of $A$ can be
viewed as $(A \# KG^*)$-$K[X]$-bimodules $M_i$ and $(A \# KG^*)$-$K
\langle x, y \rangle$-bimodule $M$ in the definitions of tameness
and wildness of $A \# KG^*$ respectively, and vice versa. Thus the
theorem follows. \end{proof}

\begin{corollary} \label{cgrungrreptype} Suppose $\ch K \nmid |G|$.
Then a finite-dimensional $G$-graded $K$-algebra $A$ is
gr-representation-finite (resp. gr-tame, gr-wild) if and only if
(forgotten $G$-grading) $A$ is representation-finite (resp. tame,
wild). \end{corollary}

\begin{proof} It follows from \cite[Theorem 4.5]{L} that, in the
case of $\ch K \nmid |G|$, $A\# KG^*$ and $A$ (forgotten
$G$-grading) have the same representation type. By
Proposition~\ref{preptypecoincide}, we are done. \end{proof}

\begin{theorem} \label{tDrozd} {\bf (Graded version of Drozd's
theorem)} A finite-dimensional $G$-graded $K$-algebra $A$ is either
graded tame or graded wild, and not both.
\end{theorem}

\begin{proof} Owing to Proposition~\ref{preptypecoincide}, it suffices
to apply Drozd's Tame-Wild Theorem (ref. \cite[Corollary C]{CB}) to
$A\# KG^*$. \end{proof}

\begin{remark} The condition $\ch K \nmid |G|$ in Corollary~\ref{cgrungrreptype}
is necessary. Indeed, let $A=KQ/\i$ be the algebra given by quiver
$$\unitlength=1mm
\begin{picture}(20,10)
\mput(0,5)(10,0){3}{\circle{2}} \mput(1,4)(10,0){2}{\vector(1,0){8}}
\mput(1,6)(10,0){2}{\vector(1,0){8}} \put(3,0){$a_3$}
\put(13,0){$a_4$} \put(3,7){$a_1$} \put(13,7){$a_2$} \put(-1,7){1}
\put(9,7){2} \put(19,7){3}
\end{picture}$$ with relations $a_1a_2=a_3a_4$ and $a_3a_2=a_1a_4$.
Then the algebra $A$ is a $\mathbb{Z}_2$-graded algebra defined by
$\deg e_i := 0$ for $i=1,2,3$, $\deg a_i := 0$ for $i=1,2$, and
$\deg a_i := 1$ for $i=3,4$. Thus the smash product $A \#
K\mathbb{Z}_2^*$ is isomorphic to the algebra $B=KQ'/\i '$ given by
quiver
$$\unitlength=1mm
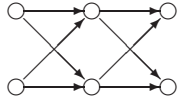
\begin{picture}(20,10)
\mput(0,0)(10,0){3}{\circle{2}} \mput(0,10)(10,0){3}{\circle{2}}
\mput(1,0)(10,0){2}{\vector(1,0){8}}
\mput(1,10)(10,0){2}{\vector(1,0){8}}
\mput(1,1)(10,0){2}{\vector(1,1){8}}
\mput(1,9)(10,0){2}{\vector(1,-1){8}}
\end{picture}$$ with all commutative relations. If $\ch K = 2$ then,
by \cite[Section 3]{GP}, $A$ is wild but $A \# K\mathbb{Z}_2^*$ is
tame.
\end{remark}

\section{Classification of acyclic superspecies}

In this section, we shall classify all acyclic superspecies
according to their graded representation type in terms of their
quivers on one hand and their superquivers on the other hand. From
now on, the underlying field $K$ is always assumed to be
algebraically closed.

\subsection{Quiver of a superspecies}

\indent\indent By Proposition~\ref{psuperdivision}, we know a
finite-dimensional gr-division superalgebra is either $K$ or $D$. In
$K$, we define $\varepsilon^K_{\bar{0}} := 1$ and
$\varepsilon^K_{\bar{1}} := 0$. In $D$, we define
$\varepsilon^D_{\bar{0}} := 1$ and $\varepsilon^D_{\bar{1}} :=
\varepsilon$.

\begin{proposition} \label{pdegree0freebasis} Each graded $D$-module has a
degree $\bar{0}$ gr-free basis.
\end{proposition}

\begin{proof} It follows from \cite[Proposition 4.6.1]{NV} that each $D$-module $M$
has a homogeneous gr-free basis $\{m_{l_M} | l_M \in I_M \}$. If
$\deg m_{l_M} = z_{l_M} \in \mathbb{Z}_2$ then $M = \oplus _{l_M \in
I_M} Dm_{l_M} = \oplus _{l_M \in I_M} D
\varepsilon^D_{-z_{l_M}}m_{l_M}$. Thus
$\{\varepsilon^D_{-z_{l_M}}m_{l_M} | l_M \in I_M \}$ is a degree
$\bar{0}$ gr-free basis of $M$. \end{proof}

Let $A$ be a finite-dimensional gr-division superalgebra. By Axiom
of Choice, we may fix a homogeneous gr-free basis $\{m_{l_M} | l_M
\in I_M \}$ for each finitely generated graded $A$-module $M$. It
follows from Proposition~\ref{pdegree0freebasis} that, in the case
of $A=D$, we may fix a degree $\bar{0}$ gr-free basis for each
graded $D$-module.

Let $\mathscr{S} = (D_i, \, _jM_i)_{i, j \in I}$ be a superspecies.
By Proposition~\ref{psuperdivision}, we know $D_i = K$ or $D$ for
all $i \in I$. Let $\{m_{l_{\, _jM_i}} | l_{\, _jM_i} \in I_{\,
_jM_i}\}$ be the fixed gr-free basis of the graded $D_j$-module $\,
_jM_i$ for all $i,j \in I$, and $z_{l_{\, _jM_i}} := \deg m_{l_{\,
_jM_i}}$. Let $I^z_{\, _jM_i} := \{l_{\, _jM_i} \in I_{\, _jM_i} |
z_{l_{\, _jM_i}} = z \}$ for all $i,j \in I$ and $z \in
\mathbb{Z}_2$. Note that, in the case of $D_j=D$, we have
$I^{\bar{0}}_{\, _jM_i} = I_{\, _jM_i}$ and $I^{\bar{1}}_{\, _jM_i}
= \emptyset$.

\begin{definition} The \mbox{\bf quiver $Q_\mathscr{S}$ of a superspecies
$\mathscr{S}$}
is the quiver $Q=(Q_0,Q_1)$ defined as follows: The vertex set $Q_0
:= \{(i,z)| i \in I, z \in \mathbb{Z}_2 \}$. We put $\{(i,
\bar{1})\} = \emptyset$ in the case of $D_i = D$. The arrow set $Q_1
:= \{a^z_{l_{\, _jM_i}} | i,j \in I, l_{\, _jM_i} \in I_{\, _jM_i},
z \in \mathbb{Z}_2 \}$ where $a^z_{l_{\, _jM_i}} : (i,z) \rightarrow
(j, z+z_{l_{\, _jM_i}})$ except for $a^{\bar{1}}_{l_{\, _jM_i}} :
(i,\bar{1}) \rightarrow (j,\bar{0})$ in the case of $D_i=K$, $D_j=D$
and $l_{\, _jM_i} \in I_{\, _jM_i}$. \end{definition}

\begin{theorem} \label{tspeciesquiver} The categories $\rep \mathscr{S}$ and $\rep Q_\mathscr{S}$ are
equivalent. \end{theorem}

\begin{proof} First of all, we define a functor ${\cal H} : \rep
\mathscr{S} \rightarrow \rep Q_{\s}$.

Let $(V_i, \, _j\phi_i)_{i,j \in I}$ be a finite-dimensional
representation of $\s$, $\{v_{l_{V_i}} | l_{V_i} \in I_{V_i}\}$ the
fixed homogeneous gr-free basis of graded $D_i$-module $V_i$ for all
$i \in I$, and $z_{l_{V_i}} := \deg v_{l_{V_i}}$. Let $I^z_{V_i} :=
\{l_{V_i} \in I_{V_i} | z_{l_{V_i}} = z \}$ for all $i,j \in I$ and
$z \in \mathbb{Z}_2$. Note that $I^{\bar{0}}_{V_i} = I_{V_i}$ and
$I^{\bar{1}}_{V_i} = \emptyset$ in the case of $D_i=D$. Since $\,
_jM_i \otimes _{D_i} V_i = \bigoplus_{l_{\, _jM_i} \in I_{\, _jM_i}}
\bigoplus_{l_{V_i} \in I_{V_i}}D_j \, m_{l_{\, _jM_i}} \otimes
v_{l_{V_i}}$ and $V_j = \bigoplus_{l_{V_j} \in I_{V_j}} D_j \,
v_{l_{V_j}}$, we have $_j\phi_i = (_j\phi_i^{l_{\,
_jM_i}l_{V_i}l_{V_j}})_{l_{\, _jM_i},l_{V_i},l_{V_j}}$ where
$_j\phi_i^{l_{_jM_i}l_{V_i}l_{V_j}} = \lambda_{l_{\, _jM_i}l_{V_i}}
\, _j\phi_i \,\rho_{l_{V_j}}$ with $\lambda_{l_{\, _jM_i}l_{V_i}}$
and $\rho_{l_{V_j}}$ the natural injection and projection. Since
$_j\phi_i$ is a graded $D_j$-module map, we have $_j\phi_i^{l_{\,
_jM_i}l_{V_i}l_{V_j}}(m_{l_{\, _jM_i}} \otimes v_{l_{V_i}}) =
c^{l_{\, _jM_i}}_{l_{V_i}l_{V_j}} \varepsilon^{D_j}_{z_{l_{\,
_jM_i}}+z_{l_{V_i}}-z_{l_{V_j}}}v_{l_{V_j}}$ where $c^{l_{\,
_jM_i}}_{l_{V_i}l_{V_j}} \in K$. We put
$c^{l_{_jM_i}}_{l_{V_i}l_{V_j}} = 0$ in the case of
$\varepsilon^{D_j}_{z_{l_{_jM_i}}+z_{l_{V_i}}-z_{l_{V_j}}}=0$. We
define the representation ${\cal H}(V_i, \, _j\phi_i)$ by ${\cal
H}(V_i, \, _j\phi_i)_{(i,z)} := K^{|I^z_{V_i}|}$ for all $(i,z)$ in
$Q_0$ and ${\cal H}(V_i, \, _j\phi_i)_{a^z_{l_{\, _jM_i}}} :=
(c^{l_{\, _jM_i}}_{l^z_{V_i}l^{z'}_{V_j}}) _{l^z_{V_i},l^{z'}_{V_j}}
: K^{|I^z_{V_i}|} \rightarrow K^{|I^{z'}_{V_j}|}$ for all
$a^z_{l_{\, _jM_i}} : (i,z) \rightarrow (j, z')$ in $Q_1$.

Let $(\az_i)_{i \in I}: (V_i, \, _j\phi_i)_{i,j \in I} \rightarrow
(W_i, \, _j\psi_i)_{i,j \in I}$ be a morphism. Then we have
commutative diagrams
$$\begin{array}{rcl} _jM_i
\otimes_{D_i} V_i & \stackrel{_j\phi_i}{\longrightarrow} & V_j\\
\downarrow 1 \otimes \az_i & & \downarrow \az_j \\ _jM_i
\otimes_{D_i} W_i & \stackrel{_j\psi_i}{\longrightarrow} & W_j
\end{array}$$ for all $i,j \in I$. Let $\az_i^{l_{V_i}l_{W_i}} := \lambda_{l_{V_i}} \, \az_i \,
\rho_{l_{W_i}}: D_iv_{l_{V_i}} \rightarrow D_iw_{l_{W_i}}$ with
$\lambda_{l_{V_i}}$ and $\rho_{l_{W_i}}$ the natural injection and
projection. Since $\az_i$ is a graded $D_i$-module map,
$\az_i^{l_{V_i}l_{W_i}}(v_{l_{V_i}}) = c_{l_{V_i}l_{W_i}}
\varepsilon^{D_i}_{z_{l_{V_i}}-z_{l_{W_i}}}w_{l_{W_i}}$ where
$c_{l_{V_i}l_{W_i}} \in K$. We put $c_{l_{V_i}l_{W_i}} = 0$ in the
case of $\varepsilon^{D_i}_{z_{l_{V_i}}-z_{l_{W_i}}}=0$. Thus we
have commutative diagrams
$$\begin{array}{rcl} \bigoplus_{l_{\, _jM_i}}
\bigoplus_{l_{V_i}}  D_j \, m_{l_{\, _jM_i}} \otimes v_{l_{V_i}} &
\stackrel{(_j\phi_i^{l_{\, _jM_i}l_{V_i}l_{V_j}}) _{l_{\,
_jM_i},l_{V_i},l_{V_j}}}{\longrightarrow} &
\bigoplus_{l_{V_j}} D_j \, v_{l_{V_j}}\\
\downarrow (1 \otimes \az_i^{l_{V_i}l_{W_i}})_{l_{V_i}, l_{W_i}} & &
\downarrow (\az_j^{l_{V_j}l_{W_j}})_{l_{V_j}, l_{W_j}} \\
\bigoplus_{l_{\, _jM_i}} \bigoplus_{l_{W_i}} D_j \, m_{l_{\, _jM_i}}
\otimes w_{l_{W_i}} & \stackrel{(_j\psi_i^{l_{\, _jM_i}l_{W_i}
l_{W_j}}) _{l_{_jM_i},l_{W_i},l_{W_j}}}{\longrightarrow} &
\bigoplus_{l_{W_j}} D_j \, w_{l_{W_j}}
\end{array}$$ for all $i,j \in I$.
We define ${\cal H}((\az_i)_{i \in I})_{(i,z)} :=
(c_{l^z_{V_i}l^z_{W_i}})_{l^z_{V_i},l^z_{W_i}}$ for all $(i,z) \in
Q_0$.

Now we show that ${\cal H}((\az_i)_{i \in I}) :=({\cal H}((\az_i)_{i
\in I})_{(i,z)})_{(i,z) \in Q_0}$ is a morphism from representation
${\cal H}((V_i, \, _j\phi_i)_{i,j \in I})$ to representation ${\cal
H}((W_i, \, _j\psi_i)_{i,j \in I})$.

Since $1 \otimes \az_i$ is diagonal on $l_{_jM_i}$, we have
commutative diagrams
$$\begin{array}{rcl}
\bigoplus_{l_{V_i}}  D_j \, m_{l_{\, _jM_i}} \otimes v_{l_{V_i}} &
\stackrel{(_j\phi_i^{l_{\, _jM_i}l_{V_i}l_{V_j}})
_{l_{V_i},l_{V_j}}}{\longrightarrow} &
\bigoplus_{l_{V_j}} D_j \, v_{l_{V_j}}\\
\downarrow (1 \otimes \az_i^{l_{V_i}l_{W_i}})_{l_{V_i}, l_{W_i}} & &
\downarrow (\az_j^{l_{V_j}l_{W_j}})_{l_{V_j}, l_{W_j}} \\
\bigoplus_{l_{W_i}} D_j \, m_{l_{\, _jM_i}} \otimes w_{l_{W_i}} &
\stackrel{(_j\psi_i^{l_{\, _jM_i}l_{W_i}l_{W_j}})
_{l_{W_i},l_{W_j}}}{\longrightarrow} & \bigoplus_{l_{W_j}}
D_jw_{l_{W_j}}
\end{array}$$ for all $i,j \in I$ and $l_{_jM_i} \in I_{_jM_i}$.

Furthermore, we have commutative diagrams
$$\begin{array}{lcl}
D_j^{|I_{V_i}|} & \stackrel{(c^{l_{\, _jM_i}}_{l_{V_i}l_{V_j}}
\varepsilon^{D_j}_{z_{l_{_jM_i}}+z_{l_{V_i}}-z_{l_{V_j}}})
_{l_{V_i},l_{V_j}}}{\longrightarrow} & D_j^{|I_{V_j}|}\\
\downarrow (c_{l_{V_i}l_{W_i}}
\varepsilon^{D_i}_{z_{l_{V_i}}-z_{l_{W_i}}}) _{l_{V_i}, l_{W_i}} & &
\downarrow (c_{l_{V_j}l_{W_j}}
\varepsilon^{D_j}_{z_{l_{V_j}}-z_{l_{W_j}}}) _{l_{V_j}, l_{W_j}} \\
D_j^{|I_{W_i}|} & \stackrel{(c^{l_{\, _jM_i}}_{l_{W_i}l_{W_j}}
\varepsilon^{D_j}_{z_{l_{\, _jM_i}}+z_{l_{W_i}}-z_{l_{W_j}}})
_{l_{W_i},l_{W_j}}}{\longrightarrow} & D_j^{|I_{W_j}|}
\end{array}$$ denoted by $(*)$, for all $i,j \in I$ and $l_{\, _jM_i} \in I_{\, _jM_i}$.

\medskip

{\bf Case $D_i=D_j=D$:} Since all $m_{l_{_jM_i}}$, $v_{l_{V_i}}$,
$v_{l_{V_j}}$, $w_{l_{W_i}}$ and $w_{l_{W_j}}$ are of degree
$\bar{0}$, we know all $\varepsilon^{D_j}_{z_{l_{\,
_jM_i}}+z_{l_{V_i}}-z_{l_{V_j}}}$, $\varepsilon^{D_j}_{z_{l_{\,
_jM_i}}+z_{l_{W_i}}-z_{l_{W_j}}}$,
$\varepsilon^{D_i}_{z_{l_{V_i}}-z_{l_{W_i}}}$ and
$\varepsilon^{D_j}_{z_{l_{V_j}}-z_{l_{W_j}}}$ in $(*)$ are equal to
1. From $(*)$ we obtain commutative diagram
$$\begin{array}{lcl}
D^{|I_{V_i}|} & \stackrel{(c^{l_{\, _jM_i}}_{l_{V_i}l_{V_j}})
_{l_{V_i},l_{V_j}}}{\longrightarrow} & D^{|I_{V_j}|}\\
\downarrow (c_{l_{V_i}l_{W_i}}) _{l_{V_i}, l_{W_i}} & & \downarrow
(c_{l_{V_j}l_{W_j}}) _{l_{V_j}, l_{W_j}} \\
D^{|I_{W_i}|} & \stackrel{(c^{l_{\, _jM_i}}_{l_{W_i}l_{W_j}})
_{l_{W_i},l_{W_j}}}{\longrightarrow} & D^{|I_{W_j}|}
\end{array}$$ Furthermore, we have commutative diagram
$$\begin{array}{lcl}
K^{|I_{V_i}|} & \stackrel{(c^{l_{\, _jM_i}}_{l_{V_i}l_{V_j}})
_{l_{V_i},l_{V_j}}}{\longrightarrow} & K^{|I_{V_j}|}\\
\downarrow (c_{l_{V_i}l_{W_i}}) _{l_{V_i}, l_{W_i}} & & \downarrow
(c_{l_{V_j}l_{W_j}}) _{l_{V_j}, l_{W_j}} \\
K^{|I_{W_i}|} & \stackrel{(c^{l_{\, _jM_i}}_{l_{W_i}l_{W_j}})
_{l_{W_i},l_{W_j}}}{\longrightarrow} & K^{|I_{W_j}|}
\end{array}$$

\medskip

{\bf Case $D_i=D$ and $D_j=K$:} Since all $v_{l_{V_i}}$ and
$w_{l_{W_i}}$ are of degree $\bar{0}$, we know all
$\varepsilon^{D_i}_{z_{l_{V_i}}-z_{l_{W_i}}}$ in $(*)$ are equal to
1. Moreover, $\varepsilon^{D_i}_{z_{l_{\,
_jM_i}}+z_{l_{V_i}}-z_{l_{V_j}}} = 1$ (resp.
$\varepsilon^{D_i}_{z_{l_{\, _jM_i}}+z_{l_{W_i}}-z_{l_{W_j}}} = 1$,
$\varepsilon^{D_i}_{z_{l_{V_j}}-z_{l_{W_j}}}=1$) precisely when
$z_{l_{\, _jM_i}} = z_{l_{V_j}}$ (resp. $z_{l_{\, _jM_i}} =
z_{l_{W_j}}$, $z_{l_{V_j}} = z_{l_{W_j}}$), and 0 otherwise. From
$(*)$ we obtain commutative diagrams
$$\begin{array}{lcl}
K^{|I_{V_i}|} & \stackrel{(c^{l^z_{\, _jM_i}}_{l_{V_i}l^z_{V_j}})
_{l_{V_i},l^z_{V_j}}}{\longrightarrow} & K^{|I^z_{V_j}|}\\
\downarrow (c_{l_{V_i}l_{W_i}}) _{l_{V_i},l_{W_i}} & & \downarrow
(c_{l^z_{V_j}l^z_{W_j}}) _{l^z_{V_j}, l^z_{W_j}} \\
K^{|I_{W_i}|} & \stackrel{(c^{l^z_{\, _jM_i}}_{l_{W_i}l^z_{W_j}})
_{l_{W_i},l^z_{W_j}}}{\longrightarrow} & K^{|I^z_{W_j}|}
\end{array}$$ for all $z \in \mathbb{Z}_2$.

\medskip

{\bf Case $D_i=K$ and $D_j=D$:} Since all $m_{l_{\, _jM_i}},
v_{l_{V_j}}$ and $w_{l_{W_j}}$ are of degree $\bar{0}$, we know all
$\varepsilon^{D_j}_{z_{l_{V_j}}-z_{l_{W_j}}}$ in $(*)$ are equal to
1. Moreover, $\varepsilon^{D_j}_{z_{l_{V_i}}-z_{l_{W_i}}}=1$ if
$z_{l_{V_i}} = z_{l_{W_i}}$, and 0 otherwise. We also have
$\varepsilon^{D_j}_{z_{l_{_jM_i}}+z_{l_{V_i}}-z_{l_{V_j}}} = 1$
(resp. $\varepsilon^{D_j}_{z_{l_{_jM_i}}+z_{l_{W_i}}-z_{l_{W_j}}} =
1$) precisely when $z_{l_{V_i}}=\bar{0}$ (resp.
$z_{l_{W_i}}=\bar{0}$), and $\varepsilon$ otherwise. From $(*)$ we
obtain commutative diagrams
$$\begin{array}{lcl}
D^{|I^z_{V_i}|} & \stackrel{(c^{l_{\, _jM_i}}_{l^z_{V_i}l_{V_j}}
\varepsilon^D_z)
_{l^z_{V_i},l_{V_j}}}{\longrightarrow} & D^{|I_{V_j}|}\\
\downarrow (c_{l^z_{V_i}l^z_{W_i}}) _{l^z_{V_i}, l^z_{W_i}} & &
\downarrow
(c_{l_{V_j}l_{W_j}}) _{l_{V_j}, l_{W_j}} \\
D^{|I^z_{W_i}|} & \stackrel{(c^{l_{\, _jM_i}}_{l^z_{W_i}l_{W_j}}
\varepsilon^D_z) _{l^z_{W_i},l_{W_j}}}{\longrightarrow} &
D^{|I_{W_j}|}
\end{array}$$ for all $z \in \mathbb{Z}_2$. Furthermore, we have commutative
diagrams
$$\begin{array}{lcl}
K^{|I^z_{V_i}|} & \stackrel{(c^{l_{\, _jM_i}}_{l^z_{V_i}l_{V_j}})
_{l^z_{V_i},l_{V_j}}}{\longrightarrow} & K^{|I_{V_j}|}\\
\downarrow (c_{l^z_{V_i}l^z_{W_i}}) _{l^z_{V_i}, l^z_{W_i}} & &
\downarrow
(c_{l_{V_j}l_{W_j}}) _{l_{V_j}, l_{W_j}} \\
K^{|I^z_{W_i}|} & \stackrel{(c^{l_{\, _jM_i}}_{l^z_{W_i}l_{W_j}})
_{l^z_{W_i},l_{W_j}}}{\longrightarrow} & K^{|I_{W_j}|}
\end{array}$$ for all $z \in \mathbb{Z}_2$.

\medskip

{\bf Case $D_i=D_j=K$:} We know
$\varepsilon^{D_i}_{z_{l_{V_i}}-z_{l_{W_i}}}=1$ (resp.
$\varepsilon^{D_j}_{z_{l_{\, _jM_i}}+z_{l_{V_i}}-z_{l_{V_j}}} = 1$,
$\varepsilon^{D_j}_{z_{l_{\, _jM_i}}+z_{l_{W_i}}-z_{l_{W_j}}} = 1$,
$\varepsilon^{D_j}_{z_{l_{V_j}}-z_{l_{W_j}}}=1$) in $(*)$ precisely
when $z_{l_{V_i}} = z_{l_{W_i}}$ (resp. $z_{l_{\,
_jM_i}}+z_{l_{V_i}} = z_{l_{V_j}}$, $z_{l_{\,
_jM_i}}+z_{l_{W_i}}=z_{l_{W_j}}$, $z_{l_{V_j}} = z_{l_{W_j}}$), and
0 otherwise. From $(*)$ we obtain commutative diagrams
$$\begin{array}{lcl}
K^{|I^z_{V_i}|} & \stackrel{(c^{l^{z'}_{\,
_jM_i}}_{l^z_{V_i}l^{z+z'}_{V_j}})
_{l^z_{V_i},l^{z+z'}_{V_j}}}{\longrightarrow} & K^{|I^{z+z'}_{V_j}|}\\
\downarrow (c_{l^z_{V_i}l^z_{W_i}}) _{l^z_{V_i}, l^z_{W_i}} & &
\downarrow
(c_{l^{z+z'}_{V_j}l^{z+z'}_{W_j}}) _{l^{z+z'}_{V_j},l^{z+z'}_{W_j}} \\
K^{|I^z_{W_i}|} & \stackrel{(c^{l^{z'}_{\,
_jM_i}}_{l^z_{W_i}l^{z+z'}_{W_j}})
_{l^z_{W_i},l^{z+z'}_{W_j}}}{\longrightarrow} & K^{|I^{z+z'}_{W_j}|}
\end{array}$$ for all $z,z' \in \mathbb{Z}_2$.

\medskip

By the above argument, we always have commutative diagram
$$\begin{array}{lcl}
K^{|I^z_{V_i}|} & \stackrel{(c^{l_{\,
_jM_i}}_{l^z_{V_i}l^{z'}_{V_j}})
_{l^z_{V_i},l^{z'}_{V_j}}}{\longrightarrow} & K^{|I^{z'}_{V_j}|}\\
\downarrow (c_{l^z_{V_i}l^z_{W_i}}) _{l^z_{V_i}, l^z_{W_i}} & &
\downarrow
(c_{l^{z'}_{V_j}l^{z'}_{W_j}}) _{l^{z'}_{V_j},l^{z'}_{W_j}} \\
K^{|I^z_{W_i}|} & \stackrel{(c^{l_{\,
_jM_i}}_{l^z_{W_i}l^{z'}_{W_j}})
_{l^z_{W_i},l^{z'}_{W_j}}}{\longrightarrow} & K^{|I^{z'}_{W_j}|}
\end{array}$$ for each arrow $a^z_{l_{\, _jM_i}} : (i, z) \rightarrow (j, z')$ in $Q_1$.
Therefore, ${\cal H}((\az_i)_{i \in I}) := ({\cal H}((\az_i)_{i \in
I})_{(i,z)})_{(i,z) \in Q_0}$ is a morphism from the representation
${\cal H}((V_i, \, _j\phi_i)_{i,j \in I})$ to the representation
${\cal H}((W_i, \, _j\psi_i)_{i,j \in I})$.

Clearly, the functor ${\cal H}$ is fully faithful. It is also dense:
Indeed, for any representation $(V_{(i,z)}, V_{a^z_{l_{_jM_i}}})$ of
the quiver $Q_\mathscr{S}$, we may define a representation $(V_i, \,
_j\phi_i)_{i,j \in I}$ of $\s$ by $V_i:=V_{(i,\bar{0})} \oplus
V_{(i,\bar{1})}$ and $_j\phi_i(m_{l_{\, _jM_i}} \otimes
v_{(i,z)}):=V_{a^z_{l_{_jM_i}}}(v_{(i,z)})$. Obviously, ${\cal
H}((V_i, \, _j\phi_i)_{i,j \in I}) \cong (V_{(i,z)},
V_{a^z_{l_{_jM_i}}})$. Thus the functor ${\cal H}$ is an
equivalence.
\end{proof}




\subsection{Classification of acyclic superspecies, I}

\indent\indent An acyclic superspecies $\mathscr{S}$ is said to be
{\it gr-representation-finite} (resp. {\it gr-tame}, {\it gr-wild})
if $T(\mathscr{S})$ is. The following result provides a
classification of acyclic superspecies in terms of their quivers:

\begin{theorem} \label{tclassification1} An acyclic superspecies $\mathscr{S}$ is
gr-representation-finite (resp. gr-tame) if and only if its quiver
$Q_\mathscr{S}$ is Dynkin (resp. extended Dykin).
\end{theorem}

\begin{proof} It follows from \cite[Theorem
2.2]{CM} that $\mod T(\mathscr{S}) \# KG^*$ and $\gr T(\mathscr{S})$
are equivalent. By Proposition~\ref{pspeciestensor} we have $\gr
T(\mathscr{S})$ and $\rep \mathscr{S}$ are equivalent. According to
Theorem~\ref{tspeciesquiver}, $\rep \mathscr{S}$ and $\rep
Q_\mathscr{S}$ are equivalent. It follows from \cite[Theorem
1.5]{ARS} that $\rep Q_\mathscr{S}$ and $\mod KQ_\mathscr{S}$ are
equivalent. Thus the representation type of $T(\mathscr{S}) \# KG^*$
and $KQ_\mathscr{S}$ coincide. By Proposition~\ref{preptypecoincide}
and the well-known results on the representation types of quivers
(ref. \cite{Ga1,DF,N}), we are done. \end{proof}

Let $\mathscr{S}$ be an acyclic superspecies. Then its quiver
$Q_\mathscr{S}$ corresponds to a unique generalized Cartan matrix,
which corresponds to a unique Kac-Moody algebra
$\mathfrak{g}_{\mathscr{S}}$ (ref. \cite{K,K2}). By \cite[Theorem
3]{K} and Theorem~\ref{tspeciesquiver}, we have the following
result:

\begin{theorem} \label{tKactheorem} {\bf (Super version of Kac's Theorem)}
Let $\mathscr{S}$ be an acyclic superspecies. Then an indecomposable
representations of $\s$ corresponds to a positive root of the
Kac-Moody algebra $\mathfrak{g}_{\mathscr{S}}$. Moreover, a positive
real root of $\mathfrak{g}_{\mathscr{S}}$ corresponds to a unique
indecomposable representation of $\mathscr{S}$. A positive imaginary
root of $\mathfrak{g}_{\mathscr{S}}$ corresponds to infinitely many
indecomposable representations of $\mathscr{S}$.
\end{theorem}

\subsection{Classification of acyclic superspecies, II}

\begin{definition} The \mbox{\bf superquiver $Q(\s)$ of a
superspecies} $\s = $ \linebreak $(D_i, \, _jM_i)_{i, j \in I}$ is a
superquiver determined by $Q(\s)_{0K}=\{i \in I | D_i=K\},$
\linebreak $Q(\s)_{0D}=\{i \in I | D_i=D\},
Q(\s)_{1\bar{0}}=\{a_{l^{\bar{0}}_{_jM_i}}: i \rightarrow j |
l^{\bar{0}}_{_jM_i} \in I^{\bar{0}}_{_jM_i} \}$ and
$Q(\s)_{1\bar{1}}=\{a_{l^{\bar{1}}_{_jM_i}}: i \rightarrow j |
l^{\bar{1}}_{_jM_i} \in I^{\bar{1}}_{_jM_i} \}$ whose elements
correspond to white vertices, black vertices, solid arrows and
dotted arrows respectively. \end{definition}

\begin{remark} Note that not every superquiver is the one of some
superspecies: For example, the superquiver $\unitlength=1mm
\begin{picture}(11,2) \put(1,2){\circle*{2}}
\put(1,2){\vector(1,0){8}} \put(10,2){\circle{2}}
\end{picture}$ cannot be the superquiver of any
superspecies.
\end{remark}

The following result provides the classification of all acyclic
superspecies according to their representation types in terms of
their superquivers:

\begin{theorem}\label{tclassification2} Let $\s=(D_i, \, _jM_i)_{i, j \in I}$ be
an acyclic superspecies over an algebraically closed field $K$. Then

(1) $\s$ is gr-representation-finite if and only if its superquiver
is the disjoint union of some superquivers of the following types:

$$\unitlength=1mm \begin{picture}(120,65)

\put(0,60){1-color superquiver whose underlying diagram is $A_n,
D_n, E_6, E_7$ or $E_8$}

\put(0,50){A(n,1):} \put(20,50){\circle{2}}
\put(21,50){\line(1,0){8}} \put(30,50){\circle{2}}
\put(33,50){$\ldots$} \put(40,50){\circle{2}}
\put(41,50){\line(1,0){8}} \put(50,50){\circle{2}}
\put(51,50){\vector(1,0){8}} \put(60,50){\circle*{2}}

\put(0,40){A(1,n):} \put(20,40){\circle{2}}
\put(21,40){\vector(1,0){8}} \put(30,40){\circle*{2}}
\put(31,40){\line(1,0){8}} \put(40,40){\circle*{2}}
\put(43,40){$\ldots$} \put(50,40){\circle*{2}}
\put(51,40){\line(1,0){8}} \put(60,40){\circle*{2}}

\put(0,30){A(2,2):} \put(20,30){\circle{2}}
\put(21,30){\line(1,0){8}} \put(30,30){\circle{2}}
\put(31,30){\vector(1,0){8}} \put(40,30){\circle*{2}}
\put(41,30){\line(1,0){8}} \put(50,30){\circle*{2}}

\put(0,20){B(1,n):} \put(20,20){\circle{2}}
\put(29,21){\vector(-1,0){8}}
\multiput(22,19)(1,0){8}{\circle*{0.1}}
\put(21,19){\vector(-1,0){0}} \put(30,20){\circle*{2}}
\put(31,20){\line(1,0){8}} \put(40,20){\circle*{2}}
\put(43,20){$\ldots$} \put(50,20){\circle*{2}}
\put(51,20){\line(1,0){8}} \put(60,20){\circle*{2}}

\put(0,10){C(n,1):} \put(20,10){\circle{2}}
\put(21,10){\line(1,0){8}} \put(30,10){\circle{2}}
\put(33,10){$\ldots$} \put(40,10){\circle{2}}
\put(41,10){\line(1,0){8}} \put(50,10){\circle{2}}
\put(59,11){\vector(-1,0){8}} \multiput(52,9)(1,0){8}{\circle*{0.1}}
\put(51,9){\vector(-1,0){0}} \put(60,10){\circle*{2}}

\put(0,0){F(2,2):} \put(20,0){\circle{2}} \put(21,0){\line(1,0){8}}
\put(30,0){\circle{2}} \put(39,1){\vector(-1,0){8}}
\multiput(32,-1)(1,0){8}{\circle*{0.1}}
\put(31,-1){\vector(-1,0){0}} \put(40,0){\circle*{2}}
\put(41,0){\line(1,0){8}} \put(50,0){\circle*{2}}
\end{picture}$$

(2) $\s$ is gr-tame but not gr-representation-finite if and only if
its superquiver is the disjoint union of some superquivers of the
following types:
$$\unitlength=1mm \begin{picture}(120,165)

\put(0,160){1-color acyclic superquiver whose underlying diagram is
$\tilde{A}_n, \tilde{D}_n, \tilde{E}_6, \tilde{E}_7$ or
$\tilde{E}_8$}

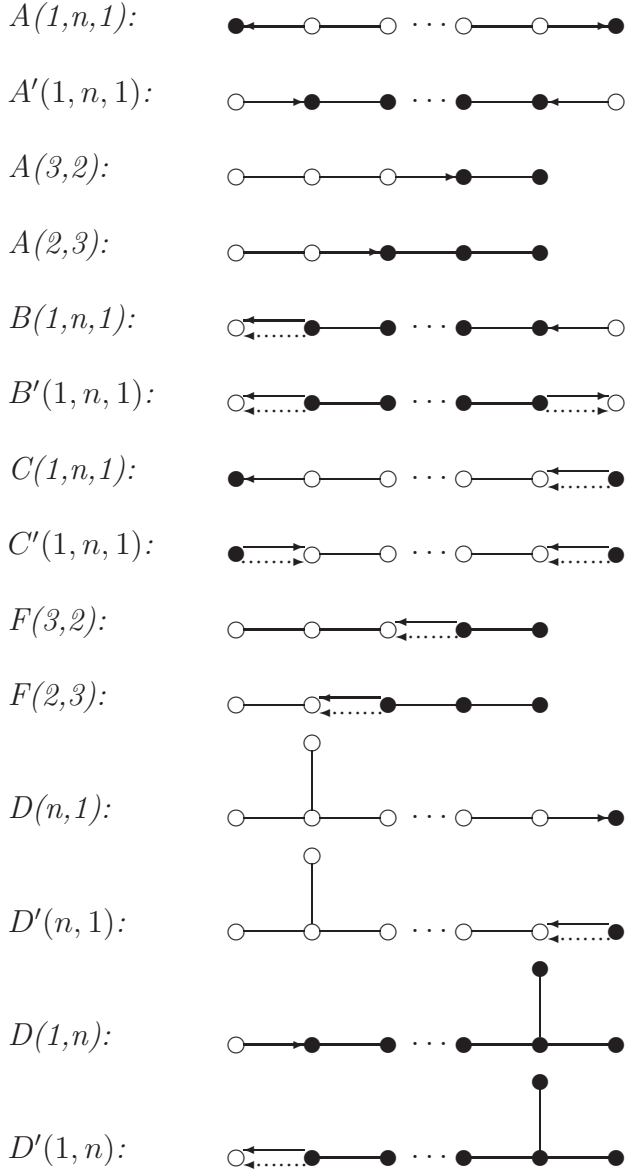
\put(0,150){A(1,n,1):} \put(30,150){\circle*{2}}
\put(39,150){\vector(-1,0){8}} \put(40,150){\circle{2}}
\put(41,150){\line(1,0){8}} \put(50,150){\circle{2}}
\put(53,150){$\ldots$} \put(60,150){\circle{2}}
\put(61,150){\line(1,0){8}} \put(70,150){\circle{2}}
\put(71,150){\vector(1,0){8}} \put(80,150){\circle*{2}}

\put(0,140){$A'(1,n,1)$:} \put(30,140){\circle{2}}
\put(31,140){\vector(1,0){8}} \put(40,140){\circle*{2}}
\put(41,140){\line(1,0){8}} \put(50,140){\circle*{2}}
\put(53,140){$\ldots$} \put(60,140){\circle*{2}}
\put(61,140){\line(1,0){8}} \put(70,140){\circle*{2}}
\put(79,140){\vector(-1,0){8}} \put(80,140){\circle{2}}

\put(0,130){A(3,2):} \multiput(30,130)(10,0){3}{\circle{2}}
\multiput(31,130)(10,0){2}{\line(1,0){8}}
\put(51,130){\vector(1,0){8}} \put(60,130){\circle*{2}}
\put(61,130){\line(1,0){8}} \put(70,130){\circle*{2}}

\put(0,120){A(2,3):} \put(30,120){\circle{2}}
\put(31,120){\line(1,0){8}} \put(40,120){\circle{2}}
\put(41,120){\vector(1,0){8}}
\multiput(50,120)(10,0){3}{\circle*{2}}
\multiput(51,120)(10,0){2}{\line(1,0){8}}

\put(0,110){B(1,n,1):} \put(30,110){\circle{2}}
\put(39,111){\vector(-1,0){8}}
\multiput(32,109)(1,0){8}{\circle*{0.1}}
\put(31,109){\vector(-1,0){0}}
\multiput(40,110)(10,0){4}{\circle*{2}} \put(41,110){\line(1,0){8}}
\put(53,110){$\ldots$} \put(61,110){\line(1,0){8}}
\put(79,110){\vector(-1,0){8}} \put(80,110){\circle{2}}

\put(0,100){$B'(1,n,1)$:} \put(30,100){\circle{2}}
\put(39,101){\vector(-1,0){8}}
\multiput(32,99)(1,0){8}{\circle*{0.1}}
\put(31,99){\vector(-1,0){0}}
\multiput(40,100)(10,0){4}{\circle*{2}} \put(41,100){\line(1,0){8}}
\put(53,100){$\ldots$} \put(61,100){\line(1,0){8}}
\put(71,101){\vector(1,0){8}}
\multiput(71,99)(1,0){8}{\circle*{0.1}} \put(79,99){\vector(1,0){0}}
\put(80,100){\circle{2}}

\put(0,90){C(1,n,1):} \put(30,90){\circle*{2}}
\put(39,90){\vector(-1,0){8}} \put(40,90){\circle{2}}
\put(41,90){\line(1,0){8}} \put(50,90){\circle{2}}
\put(53,90){$\ldots$} \put(60,90){\circle{2}}
\put(61,90){\line(1,0){8}} \put(70,90){\circle{2}}
\put(79,91){\vector(-1,0){8}}
\multiput(72,89)(1,0){8}{\circle*{0.1}}
\put(71,89){\vector(-1,0){0}} \put(80,90){\circle*{2}}

\put(0,80){$C'(1,n,1)$:} \put(30,80){\circle*{2}}
\put(31,81){\vector(1,0){8}} \multiput(31,79)(1,0){8}{\circle*{0.1}}
\put(39,79){\vector(1,0){0}} \put(40,80){\circle{2}}
\put(41,80){\line(1,0){8}} \put(50,80){\circle{2}}
\put(53,80){$\ldots$} \put(60,80){\circle{2}}
\put(61,80){\line(1,0){8}} \put(70,80){\circle{2}}
\put(79,81){\vector(-1,0){8}}
\multiput(72,79)(1,0){8}{\circle*{0.1}}
\put(71,79){\vector(-1,0){0}} \put(80,80){\circle*{2}}

\put(0,70){F(3,2):} \multiput(30,70)(10,0){3}{\circle{2}}
\multiput(31,70)(10,0){2}{\line(1,0){8}}
\put(59,71){\vector(-1,0){8}}
\multiput(52,69)(1,0){8}{\circle*{0.1}}
\put(51,69){\vector(-1,0){0}} \put(60,70){\circle*{2}}
\put(61,70){\line(1,0){8}} \put(70,70){\circle*{2}}

\put(0,60){F(2,3):} \put(30,60){\circle{2}}
\put(31,60){\line(1,0){8}} \put(40,60){\circle{2}}
\put(49,61){\vector(-1,0){8}}
\multiput(42,59)(1,0){8}{\circle*{0.1}}
\put(41,59){\vector(-1,0){0}} \multiput(50,60)(10,0){3}{\circle*{2}}
\multiput(51,60)(10,0){2}{\line(1,0){8}}

\put(0,45){D(n,1):} \put(30,45){\circle{2}}
\put(31,45){\line(1,0){8}} \put(40,45){\circle{2}}
\put(40,46){\line(0,1){8}} \put(40,55){\circle{2}}
\put(41,45){\line(1,0){8}} \put(50,45){\circle{2}}
\put(53,45){$\ldots$} \put(60,45){\circle{2}}
\put(61,45){\line(1,0){8}} \put(70,45){\circle{2}}
\put(71,45){\vector(1,0){8}} \put(80,45){\circle*{2}}

\put(0,30){$D'(n,1)$:} \put(30,30){\circle{2}}
\put(31,30){\line(1,0){8}} \put(40,30){\circle{2}}
\put(40,40){\circle{2}} \put(40,31){\line(0,1){8}}
\put(41,30){\line(1,0){8}} \put(50,30){\circle{2}}
\put(53,30){$\ldots$} \put(60,30){\circle{2}}
\put(61,30){\line(1,0){8}} \put(70,30){\circle{2}}
\put(79,31){\vector(-1,0){8}}
\multiput(72,29)(1,0){8}{\circle*{0.1}}
\put(71,29){\vector(-1,0){0}} \put(80,30){\circle*{2}}

\put(0,15){D(1,n):} \put(30,15){\circle{2}}
\put(31,15){\vector(1,0){8}} \multiput(40,15)(10,0){5}{\circle*{2}}
\put(41,15){\line(1,0){8}} \put(53,15){$\ldots$}
\multiput(61,15)(10,0){2}{\line(1,0){8}} \put(70,25){\circle*{2}}
\put(70,16){\line(0,1){8}}

\put(0,0){$D'(1,n)$:} \put(30,0){\circle{2}}
\put(39,1){\vector(-1,0){8}} \multiput(32,-1)(1,0){8}{\circle*{0.1}}
\put(31,-1){\vector(-1,0){0}} \multiput(40,0)(10,0){5}{\circle*{2}}
\put(41,0){\line(1,0){8}} \put(53,0){$\ldots$}
\multiput(61,0)(10,0){2}{\line(1,0){8}} \put(70,10){\circle*{2}}
\put(70,1){\line(0,1){8}}
\end{picture}$$
\end{theorem}

\begin{proof}
The proof of the theorem is based on Theorem~\ref{tspeciesquiver}
and the representation type classification of quivers over an
algebraically closed field (ref. \cite{Ga1, DF, N}). For some
minimal wild quivers, we refer to \cite{Ke}. Without loss of
generality, we assume that $Q(\s)$ is connected.

\textbf{Claim 1.} \textit{Let $Q(\s)$ be a 1-color superquiver. Then
$\s$ is gr-representation-finite (resp. gr-tame but not
gr-representation-finite) if and only if $\underline{Q(\s)}$ is a
Dynkin quiver (resp. extended Dynkin quiver).}

{\it Proof of Claim 1.} It follows from
Theorem~\ref{tclassification1} that, if $\underline{Q(\s)}$ is a
Dynkin quiver (resp. extended Dynkin quiver, wild quiver) then $\s$
is gr-representation-finite (resp. gr-tame but not
gr-representation-finite, gr-wild).

\textbf{Claim 2.} \textit{Let $Q(\s)$ be a 2-color superquiver with
branches, i.e., $Q(\s)$ contains vertices which have at least three
neighbors. Then $\s$ is not gr-representation-finite, and $\s$ is
gr-tame but not gr-representation-finite if and only if $Q(\s)$ is
of type $D(n, 1)$, $D'(n, 1)$, $D(1, n)$ or $D'(1, n)$.}

\textit{Proof of Claim 2.} In the case of $D_i=D$ and $D_j=K$, since
$_jM_i$ is a $\mathbb{Z}_2$-graded right $D$-module,
$m_{ij}:=\dim_{D_j}{_jM_i}$ is even. Since $Q(\s)$ has branches,
$Q(\s)$ contains a sub-superquiver $Q'$ such that $Q_0'=\{1, 2, 3,
4\}$ and $m_{i1} + m_{1i} \neq 0$ for $i=2,3,4$. If there exist $i,
j \in \{2, 3, 4\}$ such that $D_i=D_j \neq D_1$ then $Q_{\s}$
contains a subquiver of type $T_{11111}$. By
Theorem~\ref{tclassification1}, $\s$ is gr-wild. On the contrary, we
assume that either $D_1=D_2=D_3=K$ or $D_1=D_2=D_3=D$. Thus $Q(\s)$
contains a sub-superquiver $Q''$ of type $D(n, 1)$, $D'(n, 1)$,
$D(1, n)$ or $D'(1, n)$. If $Q(\s) = Q''$ then $Q_{\s}$ is of type
$\tilde{D}_{2n}$ or $\tilde{D}_{n+1}$. By
Theorem~\ref{tclassification1}, $\s$ is gr-tame but not
gr-representation-finite. If $Q''$ is a proper sub-superquiver of
$Q(\s)$ then $\s$ is gr-wild.

\textbf{Claim 3.} \textit{Let $Q(\s)$ be a 2-color superquiver
without branches and $m_{ij} \le 1$ for all $i, j \in I$. Then $\s$
is gr-representation-finite (resp. gr-tame but not
gr-representation-finite) if and only if $Q(\s)$ is of type $A(n,
1)$, $A(1, n)$, $A(2, 2)$ (resp. $A(1, n, 1)$, $A'(1, n, 1)$, $A(3,
2)$ or $A(2, 3)$).}

\textit{Proof of Claim 3.} Since $Q(\s)$ has no branches,
$\underline{Q(\s)}$ is of type $A_n$ or $\tilde{A}_n$. If
$\underline{Q(\s)}$ is of type $\tilde{A}_n$ then $Q_{\s}$ contains
a subquiver of type $\tilde{\tilde{A}}_n$. By
Theorem~\ref{tclassification1}, $\s$ is gr-wild. If
$\underline{Q(\s)}$ is of type $A_n$ then $Q(\s)$ contains a
sub-superquiver of the form $A(n, m):= \unitlength=1mm
\begin{picture}(72,2) \multiput(0,1)(10,0){4}{\circle{2}}
\multiput(1,1)(20,0){2}{\line(1,0){8}}
\multiput(13,1)(2,0){3}{\circle*{0.5}} \put(31,1){\vector(1,0){8}}
\multiput(40,1)(10,0){4}{\circle*{2}}
\multiput(41,1)(20,0){2}{\line(1,0){8}}
\multiput(53,1)(2,0){3}{\circle*{0.5}}
\end{picture}$
with $n$ white vertices and $m$ black vertices. If $Q(\s)$ is of
type $A(n, 1)$, $A(1, n)$, $A(2, 2)$ (resp. $A(1, n, 1)$, $A'(1, n,
1)$, $A(3, 2)$ or $A(2, 3)$) then $Q_{\s}$ is of type $A_{2n+1},
D_{n+2}, E_6$ (resp. $\tilde{A}_{2n+1}, \tilde{D}_{n+3},
\tilde{E}_7, \tilde{E}_6$). By Theorem~\ref{tclassification1}, $\s$
is gr-representation-finite (resp. gr-tame but not
gr-representation-finite). If $Q(\s)$ properly contains a
sub-superquiver of type $A(1, n, 1)$, $A'(1, n, 1)$, $A(3, 2)$ or
$A(2, 3)$ then $\s$ is gr-wild.

\textbf{Claim 4.} \textit{Let $Q(\s)$ be a 2-color superquiver
without branches and $m_{ij} = 2$ for some $i, j \in I$. Then $\s$
is gr-representation-finite (resp. gr-tame but not
gr-representation-finite) if and only if $Q(\s)$ is of type $B(1,
n)$, $C(n, 1)$ or $F(2, 2)$ (resp. $B(1, n, 1)$, $B'(1, n, 1)$,
$C(1, n, 1)$, $C'(1, n, 1)$, $F(3, 2)$ or $F(2, 3)$).}

\textit{Proof of Claim 4.} Note that $Q(\s)$ is a 2-color
superquiver. If $Q(\s)$ contains a sub-superquiver of type
$\unitlength=1mm \begin{picture}(12,2)
\multiput(0,1)(10,0){2}{\circle{2}}
\multiput(1,0)(0,2){2}{\line(1,0){8}}
\end{picture}$ (here a line means either a solid arrow or a dotted arrow,
certainly they are of the same orientation), $\unitlength=1mm
\begin{picture}(12,2) \put(0,1){\circle{2}}
\multiput(1,0)(0,2){2}{\vector(1,0){8}} \put(10,1){\circle*{2}}
\end{picture}$ and
$\unitlength=1mm \begin{picture}(12,2)
\multiput(0,1)(10,0){2}{\circle*{2}}
\multiput(1,0)(0,2){2}{\vector(1,0){8}}
\end{picture}$ then $Q_{\s}$ contains
a subquiver of type $\tilde{\tilde{A}}_2$. By
Theorem~\ref{tclassification1}, $\s$ is gr-wild. Suppose $Q(\s)$
contains a sub-superquiver of type $BC(n, m):= \unitlength=1mm
\begin{picture}(72,2) \multiput(0,1)(10,0){4}{\circle{2}}
\multiput(1,1)(20,0){2}{\line(1,0){8}}
\multiput(13,1)(2,0){3}{\circle*{0.5}} \put(39,2){\vector(-1,0){8}}
\put(31,0){\vector(-1,0){0}} \multiput(32,0)(1,0){8}{\circle*{0.5}}
\multiput(40,1)(10,0){4}{\circle*{2}}
\multiput(41,1)(20,0){2}{\line(1,0){8}}
\multiput(53,1)(2,0){3}{\circle*{0.5}}
\end{picture}$ with $n$ white vertices and $m$ black vertices.

If $Q(\s)$ is of type $B(1,n)$, $C(n,1)$ or $F(2, 2)$ (resp.
$B(1,n,1)$, $B'(1,n,1)$, $C(1,n,1)$, $C'(1,n,1)$, $F(3,2)$ or
$F(2,3)$) then $Q_{\s}$ is of type $D_{n+2}, A_{2n+1}$ or $E_6$
(resp. $\tilde{D}_{n+3}, \tilde{D}_{n+3}, \tilde{A}_{2n+1},
\tilde{A}_{2n+1}, \tilde{E}_7$ or $\tilde{E}_6$), and thus $\s$ is
gr-representation-finite (resp. gr-tame but not
gr-representation-finite). If $Q(\s)$ contains a proper
sub-superquiver of type $B(1,n,1)$, $B'(1,n,1)$, $C(1,n,1)$,
$C'(1,n,1)$, $F(3,2)$ or $F(2,3)$ then $\s$ is gr-wild.

\textbf{Claim 5.}  \textit{Let $Q(\s)$ be a 2-color superquiver
without branches and $m_{ij} \ge 3$ for some $i, j \in I$. Then $\s$
is gr-wild.}

\textit{Proof of Claim 5.} In this case $Q(\s)$ must contain a
sub-superquiver of type $\unitlength=1mm \begin{picture}(12,3)
\put(0,1){\circle{2}} \multiput(1,0)(0,1){3}{\line(1,0){8}}
\put(10,1){\circle{2}}
\end{picture}$, $\unitlength=1mm \begin{picture}(12,3) \put(0,1){\circle{2}}
\multiput(1,0)(0,1){3}{\vector(1,0){8}} \put(10,1){\circle*{2}}
\end{picture}$, $\unitlength=1mm \begin{picture}(12,3) \put(0,1){\circle*{2}}
\multiput(1,0)(0,1){3}{\vector(1,0){8}} \put(10,1){\circle*{2}}
\end{picture}$ or $\unitlength=1mm \begin{picture}(12,4) \put(0,1.5){\circle*{2}}
\multiput(1,2)(0,1){2}{\vector(1,0){8}}
\multiput(9,0)(0,1){2}{\vector(1,0){0}}
\multiput(1,1)(1,0){8}{\circle*{0.5}}
\multiput(1,0)(1,0){8}{\circle*{0.5}} \put(10,1.5){\circle{2}}
\end{picture}$. Thus $Q_{\s}$ contains a subquiver of type $K_3$.
By Theorem~\ref{tclassification1}, $\s$ is gr-wild.
\end{proof}

\end{document}